\documentclass[12pt]{amsart}
\usepackage[margin=1.25in]{geometry}
\usepackage[utf8]{inputenc}
\usepackage{graphicx}
\usepackage{cite}
\usepackage{enumerate}
\usepackage{xcolor}
\usepackage{setspace}
\usepackage{epstopdf}

\numberwithin{equation}{section}
\theoremstyle{plain}

\usepackage{etoolbox}
\apptocmd{\thebibliography}{\setlength{\itemsep}{10pt}}{}{}
\newtheorem{theorem}{Theorem}[section]
\newtheorem{assume}{Assumption}

\title[Long-Time Influence of Perturbations]{Long-Time Influence of Small Perturbations and Motion on the Simplex of Invariant Probability Measures}
\author[Mark Freidlin]{Mark I. Freidlin\textsuperscript{*}}\thanks{\noindent\textsuperscript{*}Dept of Mathematics, University of Maryland, College Park, MD 20742, mif@math.umd.edu}

\begin{document}
\def\e{\epsilon}
\def\Me{\mathcal{M}_{\mathrm{erg}}}
\def\M{\mathcal{M}}
\def\R{\mathbb{R}}
\def\O{\mathcal{O}}
\def\E{\mathbb{E}}
\def\d{\delta}
\def\ds{\displaystyle}
\def\sys{\{X_t\}_{t \geq 0}}
\def\hs4{ \ \ \ \ }
\newcommand{\red}[1]{\textcolor{red}{#1}}
\newcommand{\blue}[1]{\textcolor{blue}{#1}}

\maketitle
\begin{abstract}
	We present a general approach to a broad class of asymptotic problems related to the long-time influence of small perturbations, of both the deterministic and stochastic type. The main characteristic of this influence is a limiting motion on the simplex of invariant probability measures of the non-perturbed system in an appropriate time scale. We consider perturbations of dynamical systems in $\R^n$, linear and nonlinear perturbations of PDE's, wave fronts in the reaction-diffusion equations, homogenization problems and perturbations caused by small time delay. The main tools we use in these problems are limit theorems for large deviations, modified averaging principle and diffusion approximation. 
\end{abstract}

\section{Introduction}\label{sec_intro}
The goal of this paper is to present a general approach to a broad class of problems related to the long-time effects of perturbations in a variety of systems. The original system could, for instance, be a finite-dimensional dynamical system, a Markovian stochastic process or a semi-flow in a functional space related to an evolutionary PDE. The perturbations could also be of several different types: e.g. deterministic or stochastic perturbations of the equation itself or the initial condition, small delays, or perturbations of the domain in an initial-boundary problem are considered.  Many of the problems are closely related to PDE's with a small parameter.

Let the metric space $(\mathcal{E}, \rho)$ be the phase space of the original, non-perturbed system $\{X_t\}_{t \geq 0}$. We assume that $\sys$ satisfies Assumption $\ref{assume1}$, given below.
\begin{assume}\label{assume1}
	For each initial point $X_0 = x \in \mathcal{E}$, there exists an invariant measure $\mu_x$ for system $\{X_t\}_{t \geq 0}$ defined on the Borel $\sigma$-field of $(\mathcal{E},\rho)$ such that for any bounded continuous function $f:\mathcal{E} \to \R$,
\begin{equation}\label{eq_intro_ergod_prop}
	\lim_{T \to \infty} \frac{1}{T} \int_0^T f(X_t)dt = \int_{\mathcal{E}} f(y) \mu_x(dy).
\end{equation}
\end{assume}

Denote by $\M$ the simplex of all invariant probability measures of $\sys$ and let $\Me$ be the collection of all ergodic invariant probability measures. The simplex $\mathcal{M}$ is the convex envelope of $\Me$. 

Consider the following dynamical system in $\mathbb{R}^n$:
\begin{equation}\label{eq_sys_generic}
	\dot{X}_t = b(X_t), \hs4  X_0 = x \in \mathbb{R}^n.
\end{equation}
We assume that the vector field $b(x)$ is Lipschitz continuous, i.e. $|b(x)-b(y)| \leq K|x-y|$ for some $K < \infty$. Together with $\eqref{eq_sys_generic}$, consider the family of perturbed systems $\{X_t^\e\}_{t > 0}$ given by
\begin{equation}\label{eq_sys_generic_pert}
	\dot{X_t^\e} = b(X_t^\e)+\e \beta(X_t^\e), \hs4 X_0^\e = x \in \mathbb{R}^n,
\end{equation}
where $\e \ll 1$ is a small positive parameter and $\beta(x)$ is a bounded and sufficiently regular vector field. Let $N := \ds \sup_x |\beta(x)| < \infty$ and $m^\e(t):= \ds \max_{0 \leq s \leq t} |X_s^\e - X_s|$. It follows immediately from $\eqref{eq_sys_generic}$ and $\eqref{eq_sys_generic_pert}$ that 
\begin{equation*}
	m^\e(t) \leq K \int_0^t m^\e(s)ds + \e N t, 
	\end{equation*}
and hence also that
\begin{equation}\label{eq_intro_sys_bound_aftergron}
	m^\e(t) \leq \e N t e^{kt}.
\end{equation}
Inequality $\eqref{eq_intro_sys_bound_aftergron}$ implies that $\ds \lim_{\e \to 0} \ds m^\e(t(\e)) =0$ for fixed $t(\e) = t$, as well as for any $t(\e)$ tending to infinity as $\e \to 0$ which grows sufficiently slowly. Namely, if $1 \ll t(\e) \ll \alpha \log \frac{1}{\e}$ as $\e \downarrow 0$ for some $\alpha < 1 / K$, then $\ds \lim_{\e \to 0} m^\e(t(\e) = 0$.

Equation $\eqref{eq_intro_sys_bound_aftergron}$ together with $\eqref{eq_intro_ergod_prop}$ implies that $X_{t(\e)}^\e$ remains close to the support of the measure $\mu_x$. For $\e$ small and $1 \ll t(\e) \ll \e^{-\alpha}$ with $\alpha < \frac{1}{K}$, the position of $X_{t(\e)}^\e$ can be characterized by $\mu_x \in \mathcal{M}$ and a point $y^\e = y(x,t(\e))$ of the support $G_x$ of $\mu_x$ closest  to $X_{t(\e)}^\e$. 

Bounds similar to $\eqref{eq_intro_sys_bound_aftergron}$ hold if the perturbation also has a stochastic component. For instance, consider the system  
\begin{equation*}
	\dot{X}_t^\e = b(X_t^\e) + \e \beta(X_t^\e) + \sqrt{\e} \sigma(X_t^\e) \dot{W}_t^\e,\hs4  X_0^\e = x,
\end{equation*}
where $W_t$ is the Wiener process in $\R^n$ and $\sigma(x)$ is a matrix with bounded, Lipschitz continuous entries. Convergence of $m^\e(t)$ to zero in this case holds in probability. Moreover, similar bounds may hold even if the original system has a stochastic component as well.

Under certain assumptions and after an appropriate time re-scaling, we will see that the limit
\begin{equation*}
	\lim_{\e \to 0} \mu^\e_{X_{f(\e,t)}^\e} = \bar{\mu}_t^x \in \mathcal{M}.
\end{equation*}  
exists. Here $\mu^\e_x$ is the invariant probability measure of the perturbed system while $f(\e,t)$ provides the time re-scaling. For instance, in the sequel we consider cases of $f(\e,t)$ being $t/\e$, $t/\e^2$ and $\exp(t/\e)$. Note also that $\ds \lim_{\e \to 0} y(x,f(\e,t))$ exists only in special cases. However, after certain regularization, one can consider the distribution of $y(x,f(\e,t))$ and this distribution converges to the measure $\bar{\mu}_t^x$ as $\e \to 0$. Moreover, the evolution of measures $\bar{\mu}_t^x$ in time satisfies the semi-group property $\ds \bar{\mu}_{t+s}^x(\gamma)  = \int_{\mathcal{E}} \bar{\mu}_s^x(dy)\bar{\mu}_t^y(\gamma),$ and in many cases it can be described explicitly.

Our approach to studying the long-time influence of small perturbations is to some extent a generalization of the classical averaging principle. For example, consider an one-degree-of-freedom Hamiltonian system
\begin{equation}\label{eq_intro_Hamil_sys}
	\dot{X}_t = \bar{\nabla}H(X_t), \hs4 X_0 = x \in \mathbb{R}^2,
\end{equation}
where $\bar{\nabla} H = (-\frac{\partial H}{\partial x_2}, \frac{ \partial H}{\partial x_1})$ and $H(x)$ is a smooth one-well Hamiltonian. Moreover, assume that $\ds \min_{x \in \mathbb{R}^2} H(x) = H(0) = 0$, $\ds \lim_{|x| \to \infty} H(x) = \infty$ and $\nabla H(x) \neq 0$ whenever $x \neq 0$. In this case, the set $\Me$ of ergodic invariant measures can be parametrized by their value of energy $H$. In particular, on each curve $C(z) = \{x \in \mathbb{R}^2:H(x) = z\}$ there is one invariant probability measure $\mu_z$ having density
\begin{equation*}
	m_z(x) = \frac{1}{T(z)} \frac{1}{| \nabla H(x)|)}, \hs4 x \in C(z),
\end{equation*} with respect to the length element on $C(z)$. Here, $T(z) = \oint_{C(z)} \frac{d \ell}{|\nabla H(x)|}$ is the period of rotation along $C(z)$. In this sense, $z \geq 0$ can be considered as a coordinate in the set $\Me$.

Consider now the following perturbation of $\eqref{eq_intro_Hamil_sys}$:
\begin{equation}\label{eq_intro_Hamil_averaged_fast}
	\dot{X}_t^\e = \bar{\nabla} H(X_t^\e) + \e \beta(X_t^\e), \hs4 X_0^\e = x,
\end{equation}
where $\beta(x)$ is a smooth bounded vector field and $0  < \e \ll 1$ as always. The classical averaging principle (see for instance \cite{Arnold78}, Section 51) implies that under mild additional assumptions $H(X_{t/\e}^\e)$ converges as $\e \to 0$ uniformly on each finite time interval to the solution $Z(t)$ to
\begin{equation}\label{eq_intro_Hamil_averaged}
	\dot{Z}(t) = \bar{\beta}(Z(t)), \hs4 Z(0) = H(x).
\end{equation}
Here, $\bar{\beta}$ is the appropriately averaged version of $\beta$ defined by
\begin{equation*}
	\bar{\beta}(z) = \frac{1}{T(z)} \int_{G(z)} \mathrm{div} \beta(x) dx,
\end{equation*}
where $G(z)$ is the domain in $\mathbb{R}^2$ bounded by the curve $C(z)$. Equation $\eqref{eq_intro_Hamil_averaged}$ effectively determines the dynamics of $H(X_{t/\e}^\e)$ for small $\e$. Moreover, in this case $y(x,t/\e)$ changes rapidly as $\e \to 0$ and has no pointwise limit. However, one can prove that $y(x,t/\e)$ approaches, in a sense, the probability distribution on $C(z_t)$ with the density $m_{z_t}$.

Suppose now that the perturbation of $\eqref{eq_intro_Hamil_sys}$ also has a stochastic component; i.e. consider the system
\begin{equation*}
	\dot{X}_t^{\e,\d} = \bar{\nabla} H(X_t^{\e,\d}) + \e \beta(X_t^{\e,\d}) + \sqrt{\e \d} \ \sigma(X_t^{\e,\d}) \dot{W}_t, \hs4 X_0^{\e,\d} = x,
\end{equation*}
where $\d >0$ is another small parameter, $\sigma(x)$ is a $2\times 2$ matrix with smooth bounded entries and $W_t$ is the Wiener process in $\mathbb{R}^2$. In addition, set $a(x) := \sigma(x) \sigma^*(x)$, 
\begin{equation*}
 \bar{a}(z) := \frac{1}{T(z)} \int_{G(z)} \mathrm{div} (a(x) \nabla H(x)) dx,
\end{equation*}
and $\bar{\sigma}(z) := \bar{a}(z)^{1/2}$. Next, define the diffusion process $Z_t^\d$ by the equation 
\begin{equation*}
	\dot{Z}_t^\d = \bar{\beta}(Z_t^\d) + \sqrt{\d} \bar{\sigma}(Z_t^\d) \dot{W}_t, \hs4 Z_0^\d = H(x).
\end{equation*}
The point $z = 0$ is inaccessible for $Z_t^\d$ if $Z_0^\d > 0$, so that $Z_t^\d$ is a process on the set $ \{ z \geq 0 \}$, which parametrizes $\Me$. One can prove (see, for instance \cite{Freidlin12}, Ch. 8) that processes $H(X_{t/\e}^{\e,\d})$ converge weakly on each finite time interval to $Z_t^\d$ as $\e \to 0$.

Suppose now that the Hamiltonian $H(x)$ has more than one well. Then the corresponding set of ergodic invariant probability measures $\Me$ can be parametrized by the points of a graph which has interior vertices (see Section \ref{sec_graph}). It turns out that, in general, the interior vertices can be accessible in a finite time and the classical averaging principle does not work: it needs to be supplemented by a description of the behavior of the limiting trajectory on $\Me$ after hitting an interior vertex.

If a system with more than one degree-of-freedom is considered, the classical averaging principle with convergence for any fixed initial point does not hold, even in the case when global action-angle coordinates can be introduced. Due to resonances, the set $\Me$ for such a system is more sophisticated than a domain in $\mathbb{R}^n$. In this case, the convergence of the first integrals to a limiting motion in $\mathbb{R}^n$ holds only in the Lebesgue measure of the phase space \cite{Anosov60}, \cite{Lochak88}. In the classical averaging principle one can say that the phase space of the limiting process is a domain in $\R^n$. We will see that in many interesting cases, the phase space of the limiting evolution has a more sophisticated topological structure. In particular, it can be a graph or an open book space, or it can consist of just a finite number of points. Fast components of the perturbed motion roughly speaking, change on the supports of corresponding measures $\mu \in \mathcal{M}$, and, unlike the classical case, these supports can be very different. For instance, they may have different dimensions.

Throughout the remainder, we will consider deterministic and stochastic perturbations together. This is natural from the point of view of applications, but this is also very useful from a mathematical point of view. Regularization by the addition of random perturbations often makes the problem simpler. In the classical averaging principle, some regularization is used as well: namely, the replacement of the convergence starting from any fixed initial point by the convergence in the Lebesgue measure of the phase space. Studying this problem using convergence in measure is equivalent to adding a random perturbation to the initial point. Under certain conditions, this type of regularization allows for the description of the limiting evolution of the slow component in an important problem of perturbations of completely integrable Hamiltonian systems with many degrees-of-freedom in a domain where action-angle coordinates can be introduced. However, this type of regularization as a rule does not work if first integrals have critical points; for instance in the case of Hamiltonian systems with one degree-of-freedom when the Hamiltonian has more than one saddle point. Other more ``powerful'' types of stochastic regularization should be used in this case. For example, one can add random perturbations to  the equation rather than to the initial point. We will see then that the limiting motion of the slow component can be described for a broad class of problems. In particular, using such a regularization one can describe the long time behavior of system $\eqref{eq_intro_Hamil_averaged_fast}$ even if the Hamiltonian has many wells (see Section \ref{sec_graph}). It turns out that if $H(x)$ has saddle points, the long-time behavior of the perturbed system, which is purely deterministic, will be in some sense stochastic. This stochasticity is independent of the type of regularization and is an intrinsic property of the deterministic system $\eqref{eq_intro_Hamil_sys}$ and its perturbations. In fact, the stochasticity is actually caused by instabilities of the saddle points of system $\eqref{eq_intro_Hamil_sys}$.

The averaging principle is a result of the law-of-large-numbers type. In probability theory, there are, roughly speaking, three types of limit theorems: laws of large numbers, central limit theorem (or diffusion approximation), and limit theorems for large deviations. We will see that the limiting long-time behavior of the perturbed system in an appropriate time scale, should in many cases be described by a limit theorem for large deviations or by a diffusion approximation.

Let us now present our general approach for studying the long-time influences of small perturbations on a system $X_t$ with phase space $\mathcal{E}$. First, one should consider the simplex $\mathcal{M}$ of invariant probability measures of the non-perturbed system and check that Assumption \ref{assume1} is satisfied. Next, one should introduce the projection map
\begin{equation*}
	Y: \mathcal{E} \to \mathcal{M}, \hs4 Y(x) = \mu_x \in \mathcal{M}.
\end{equation*}
If $X_t(x)$ is the non-perturbed trajectory starting at $x \in \mathcal{E}$, then $Y(X_t(x)) = \mu_x$ for any $ t \geq 0$. The perturbed trajectory $X_t^\e(x)$ then induces a motion $Y_t^\e(x) = Y(X_t^\e(x))$ on $\mathcal{M}$. This motion can be rather complicated and in particular, it can have a memory. We will see that in many interesting cases there will exist a time re-scaling $f(\e,t)$ such that $Y^\e_{f(\e,t)}(x)$ converges as $\e \to 0$ to a function with values in $\mathcal{M}$. This limiting motion $\bar{\mu}_t^x$ will satisfy the semi-group property: $\bar{\mu}_{s+t}^x(\gamma) = \int_{\mathcal{E}} \bar{\mu}_s^x(dy) \bar{\mu}_t^y(\gamma)$ for each Borel set $\gamma \subset \mathcal{E}$. In certain cases, one may be able to describe $\bar{\mu}^x$ as a Markov stochastic process on the simplex $\mathcal{M}$.

My goal in this paper is to present a unified approach to various problems related to the long-time asymptotics in problems with a small parameter. Therefore, I consider often not a general case, but a simplest situation where one or another effect which we are interested in can be observed. Moreover, since complete proofs, as a rule, use specific for each problem techniques, we restrict ourselves just to sketches of proofs and list some references where the full proof can be found. There exists a large literature where many of the problems mentioned here are addressed specifically. The bibliography to this paper is incomplete, and I apologize that many interesting papers are not included in the list of references.

The outline of this paper is as follows. In Section \ref{sec_finite}, some problems where $\Me$ consists of a finite number of points are considered, and the long-time evolution is defined, roughly speaking, by transitions between asymptotically stable attractors. These transitions occur in an exponentially large time scale and are caused by large deviations. Exit problems, hierarchy of cycles (of Markov chains), metastability and related PDE problems are addressed in this section. One should mention that there are other systems with finite number of ergodic probability measures: for instance, systems without stable attractors where the time evolution goes much faster.

Systems for which $\Me$ can be parametrized by points of a graph are considered in Section \ref{sec_graph}. Regularization of the system, stochasticity of long-time behavior of pure deterministic systems and the Dirichlet problem for PDEs with a small parameter in higher derivatives are considered in this section as well as perturbations of the Landau-Lifshitz equations and area-preserving flows on the 2-torus.

In Section \ref{sec_book}, problems where the limiting motion can be described as a motion on an open book space are considered. In particular, the periodic homogenization for second order elliptic PDEs when the corresponding system on the torus has non-unique invariant probability measures is discussed.

In Section \ref{sec_flows}, we consider perturbations of semi-flows in functional spaces. Wave fronts and other patterns in Reaction-Diffusion Equations are considered. Long-time effects caused by small delay are considered as well.

Fast oscillating perturbations of finite-dimensional dynamical systems are considered in Section \ref{sec_oscil}. If the corresponding non-perturbed system has asymptotically stable attractors and the fast-oscillating noise has strong enough mixing, the transition between the attractors can occur in an exponentially long time and can be described by the corresponding action functional. But if the non-perturbed system has first integrals, the situation is different. For Hamiltonian systems with one degree of freedom and one well Hamiltonian, the limiting evolution can be described by a diffusion approximation if the noise has good enough mixing. In the case of many degrees of freedom and many first integrals, one should make additional assumptions concerning the smallness of the resonance set.

In the last Section \ref{sec_phantom}, we consider the case of the non-perturbed system for which our main Assumption \ref{assume1} is not satisfied. More precisely, we consider a system which has no finite invariant measures. It turns out that the noise itself can induce in this case invariant, in a sense, measures and attractors which define the long time behavior. This noise-induced behavior we call phantom dynamics.

\section{Systems with a finite number of ergodic invariant probability measures}\label{sec_finite}
\begin{figure}[h]

\centering
\includegraphics[scale = .9]{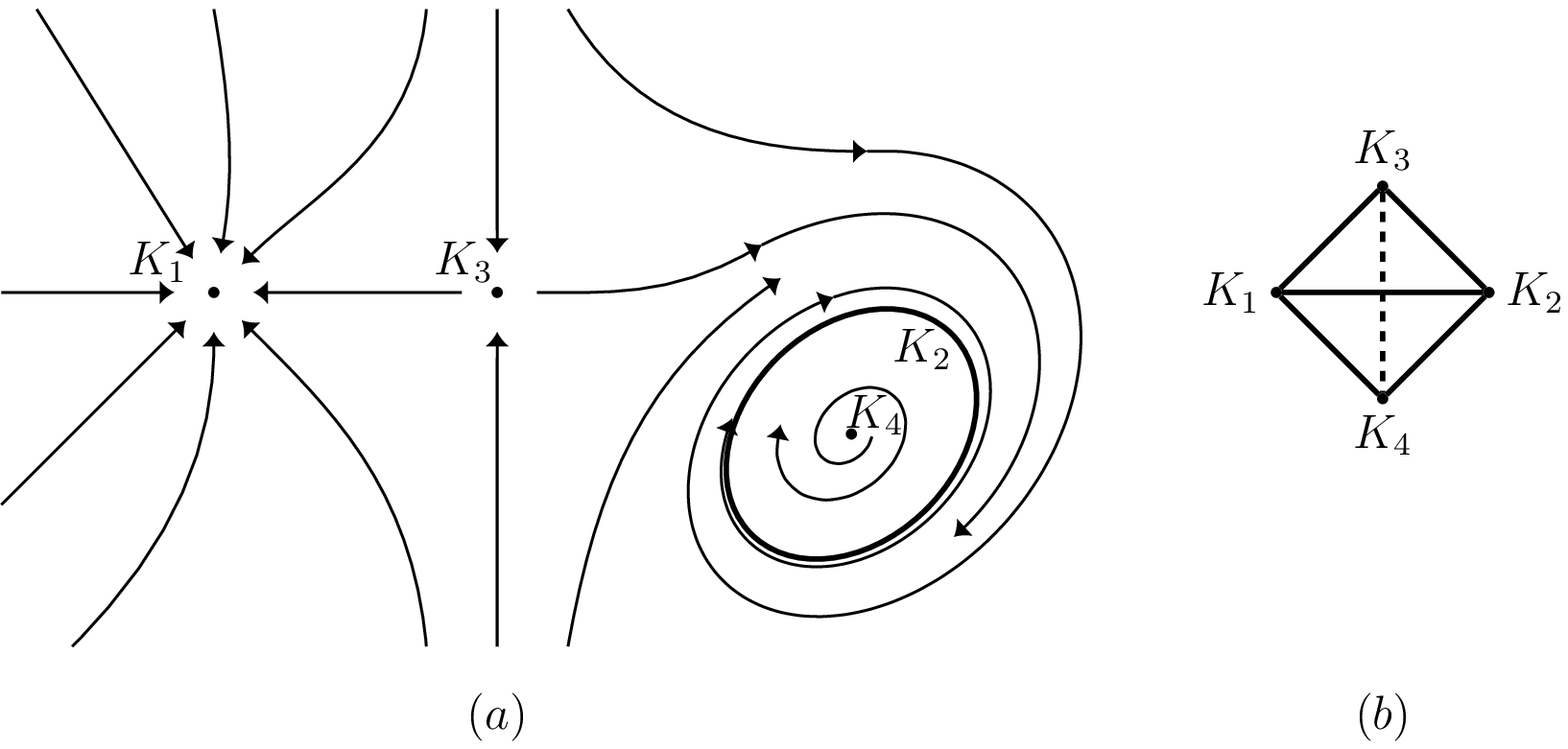}
\vspace*{-6mm}
\caption{}
\label{fig_finite_example}
\end{figure}

Consider the following dynamical system in $\mathbb{R}^n$:
\begin{equation}\label{eq_finite_dyn_sys}
	\dot{X}_t = b(X_t), \hs4 X_0 = x \in \mathbb{R}^n,
	\end{equation}
where the vector field $b:\R^n \to \R^n$ is assumed to be smooth, bounded and satisfy
\begin{equation*}
	x \cdot b(x) \leq - \alpha |x|, \hs4 \forall x : |x| \geq \tau,
\end{equation*} for some sufficiently small $\alpha > 0$ and large $\tau > 0$. This condition guarantees that all trajectories will return to the ball of radius $\tau$ quickly enough. Next, assume that there exist a finite number of invariant compact sets $K_1,...,K_{\tilde{\ell}} \subset \mathbb{R}^n$ for system $\eqref{eq_finite_dyn_sys}$ such that each trajectory is attracted as $t \to \infty$ to one of the sets $K_1,...,K_{\tilde{\ell}}$. Moreover, suppose that for some $\ell \leq \tilde{\ell}$, $K_1,...,K_{\ell}$ are asymptotically stable and that each of them supports only one ergodic invariant probability measure.  Let each of $K_i$ with $i > \ell$ belong to the union of a finite number of smooth manifolds of dimension less than $n$. A typical example is shown in Figure \ref{fig_finite_example}.

In Figure \ref{fig_finite_example}, we have here three equilibriums $K_1$, $K_3$ and $K_4$, and one invariant compact-limit cycle $K_2$; $K_1$ and $K_2$ are asymptotically stable while $K_3$ and $K_4$ are unstable equilibriums. The separatrix consists of two trajectories entering the saddle point $K_3$; all points of this separatrix are attracted to $K_3$. Points to the right of separatrix besides $K_4$ are attracted to $K_2$, points to the left are attracted to $K_1$ and the unstable point $K_3$ is attracted to itself. It is clear that Assumption \ref{assume1} is satisfied in this example, since for each $x \in \R^2$ the corresponding measure $\mu_x$ belongs to the set $ \Me$, which in this case consists of three $\d$-measures concentrated at $K_1,K_3,K_4$ and one measure at $K_2$.

Consider now the perturbed system. If the initial point $x$ is attracted to an asymptotically stable compact $K_i$ and the perturbations are the same as in equation $\eqref{eq_sys_generic_pert}$, then the perturbed trajectory will enter a neighborhood of $K_i$ and stay there forever provided that $\mathcal{\e}$ is small enough. In fact, one can show that
\begin{equation*}
	\lim_{\e \to 0} \frac{1}{t(\e)} \int_0^{t(\e)} f(X_s^\e)ds = \int_{K_i} f(y) \mu_i(dy),
\end{equation*}
for each continuous function $f:\R^n \to \R$ where the time-scale $t(\e)$ grows with $\e^{-1}$ and $\mu_i \in \Me$ is the invariant probability measure with support on $K_i$.

If the perturbations have a powerful enough stochastic component, then transitions between the asymptotically stable compacts are possible. For instance, consider the system
\begin{equation}\label{eq_finite_stoc_pert}
	\dot{X}_t^\e = b(X_t^\e)  + \e \beta(X_t^\e) + \sqrt{\e} \sigma(X_t^\e) \dot{W}_t, \hs4 X_0^\e = x \in \R^n,
\end{equation}
where $\beta$ is the same as in $\eqref{eq_sys_generic_pert}$, $W_t$ is the Wiener process in $\R^n$ and $\sigma:\R^n \to \R^{n \times n}$ has smooth bounded entries and is such that $a = \sigma \sigma^*$ is uniformly positive definite. A theory of long-time behavior of systems with asymptotically stable attractors perturbed by a Gaussian noise was developed in \cite{Freidlin12}. Here, we recall some main notions and results.

We introduce the following functional $S_{0T}(\varphi)$ on the space $C_{0T}$ of continuous functions on $[0,T]$ with values in $\mathbb{R}^n$:
\begin{equation*}
	S_{0T}(\varphi) := \begin{cases}
		\ds \frac{1}{2}  \int_0^T \Big(a^{-1}(\varphi_s)(\dot{\varphi}_s - b(\varphi_s)) \Big)\cdot \Big(\dot{\varphi}_s - b(\varphi_s)\Big) ds, & \text{if }\varphi \in AC_{0T},
		\\ +\infty, & \text{otherwise}.
	\end{cases}
\end{equation*}
Here, $AC_{0T}$ denotes the space of absolutely continuous functions on $[0,T$]. The functional $S_{0T}(\varphi)$ is called the (normalized) action functional for the family of processes $\{X_t^\e\}_{0 \leq t \leq T}$ defined by $\eqref{eq_finite_stoc_pert}$. Next, we set
\begin{equation}\label{eq_finite_quasi}
	V_{ij} := \inf \{ S_{0T}(\varphi) : \varphi \in C_{0T}, \ \varphi_0 \in K_i, \ \varphi_T \in K_j, \ T \geq 0 \}.
\end{equation}
Roughly speaking, $\exp\big( -\frac{1}{\e} S_{0T}(\varphi) \big)$ is the main term of $\mathbb{P}_{\varphi_0} (\max |X_t^\e - \varphi_t| < \d )$ for small $\d$ and $\e \downarrow 0$, while the numbers $V_{ij}$ characterize the difficulty of transition from $K_i$ to $K_j$. Rigorous statement can be found in \cite{Freidlin12}, Ch. 3, 5.

Next, define the following mapping $\mathcal{N}$ from the set $\mathcal{L} = \{1,...,\ell\}$ to itself:
\begin{equation*}
	\mathcal{N}(i) = j \text{ if }V_{ij} = \min \{V_{ik}:k \in \mathcal{L}, k \neq i \}.
\end{equation*}
In other words, the mapping $\mathcal{N}(i)$ indicates the most likely set $K_j$ to go to from $K_i$. If the minimum is achieved just for one $k = k(i)$, so that $\mathcal{N}(i)$ consists of one point, we call the system generic. In this case, suppose the initial point $x$ is attracted to the compact $K_{i_0}$. Set $i_1 =\mathcal{N}(i_0),...,i_{k+1} = \mathcal{N}(i_k),...$. The sequence $\{i_j\}$ will eventually start to repeat itself, and we will obtain a cycle $i_{m_1} \to i_{m_1+1} \to ... \to i_{m_1}$. In this way, one can get a decomposition of the set $\mathcal{L}$ into 1-cycles.
\begin{figure}[h]

\centering
\includegraphics[scale = 0.8]{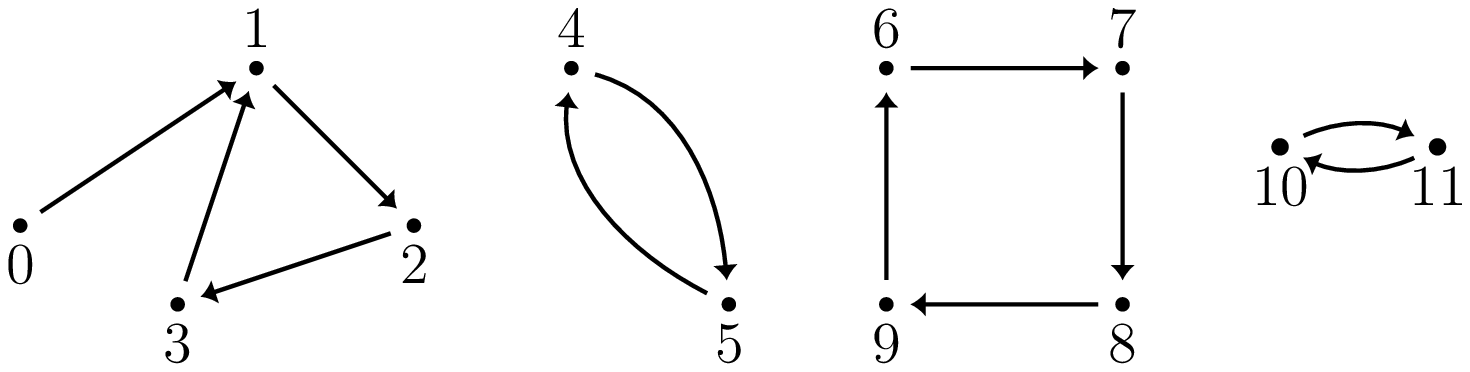}
\vspace*{-4mm}
\caption{}
\label{fig_finite_cycle}
\end{figure}

There are $\ell = 12$ asymptotically stable compacts in the system shown in Figure \ref{fig_finite_cycle}. At each point $i \in \mathcal{L}$, there is an arrow from $i \to \mathcal{N}(i)$. There are 5 1-cycles: $\{0\}, \{1 \to 2 \to 3 \to 1\}, \{4 \to 5 \to 4\}, \{6 \to 7 \to 8 \to 9 \to 6\}, \{10 \to 11 \to 10\}$. Note that if no arrow leads to a point, then such a point is considered as a separate 1-cycle. Since the transitions between the compacts $K_i$ are occurring because of the noise, the transition times are random variables. However, as it follows from the results of Ch. 6 of \cite{Freidlin12}, the logarithmic asymptotics as $\e \downarrow 0$ of these transition times are not random and can be explicitly expressed through the numbers $V_{ij}$.

Similar to construction of 1-cycles, transitions between 1-cycles can be described by 2-cycles, where each 2-cycles consists of a periodic sequence of 1-cycles. One can then define 3-cycles, and so on until all compacts $K_i$ will be involved. The structure of this hierarchy of cycles, in the generic case was described in \cite{Freidlin77} (see also Ch. 6 of \cite{Freidlin12}). It is completely defined by the numbers $V_{ij}$. These results imply the following theorem.
\begin{theorem}\label{th_finite_cycles}
	Suppose that Assumption \ref{assume1} is satisfied and let $\{\mu_i\}_{ i \in \mathcal{L}}$ be the ergodic probability measure (unique) concentrated on the compacts $\{K_i\}_{i \in \mathcal{L}}$. Assume that the system is generic and $\ds \lim_{\e \to 0} T(\e) = \infty$.	Then, for each $i \in \mathcal{L}$, there exist a set of numbers $\{\lambda_1,...,\lambda_m\}$, $0 = \lambda_0 < \lambda_1 < ... < \lambda_m$, and a set of integers $\{i_1^*,...,i_m^*\} \subset \mathcal{L}$ such that if 
	\begin{equation*}
		\lambda_\tau < \liminf_{\e \to 0} \e \log T(\e) \leq \limsup_{\e \to 0} \e \log T(\e) < \lambda_{\tau+1},
	\end{equation*}
	for some $\tau = 0,1,...,m$ (put $\lambda_{m+1} = \infty$), and the initial point $X_0^\e = x$ is attracted to $K_i$, then
	\begin{equation*}
		\lim_{\e \to 0} \frac{1}{T(\e)} \int_0^{T(\e)} f(X_s^\e) ds = \int_{K_{i_\tau^*}} f(y) \mu_{i_\tau^*}(dy).
	\end{equation*}
	The numbers $\lambda_j$ and $i_j^*$ can be expressed explicitly through the $V_{ij}$. The measure $\mu_{i_\tau^*}$ is called the metastable (or sublimit as in \cite{Freidlin77}) distribution of $X_t^\e$ for the initial point $x$ and the time scale $T(\e)$.
\end{theorem}

Applications of these types of large-deviation techniques to various physical models as well as some refinements and generalizations can be found in \cite{Olivieri05} and references therein.

If the system is not generic, for instance if $\ds \min_{k \in \mathcal{L}, k \neq i}V_{ik}$ is achieved for more than one $k$, the transition probabilities from $i$ to $k$ may have the same order as $\e \to 0$ for different values of $k$. This is called ``rough symmetry'' \cite{Freidlin14}. For systems with rough symmetry, it is not enough to use just the logarithmic asymptotics given by the action functional. One should look at the pre-exponential factor. The results concerning the pre-exponential factors can be found in \cite{Bovier15}, \cite{Landim19} and in references there.

Theorem \ref{th_finite_cycles} can also be used to obtain some results for the PDEs corresponding to the diffusions. We formulate some of them for the case of two space variables and vector field $b(x)$ as in Figure \ref{fig_finite_example}. Let $a(x) = (a_{ij}(x))$ be a positive definite $2 \times 2$ matrix with bounded twice continuously differentiable entries. Consider the following Cauchy problem on $\R^+ \times \R^2$:
\begin{equation}\label{eq_finite_Cauchy}
	\frac{\partial u^\e}{\partial t}(t,x) = \frac{\e}{2} \sum_{i,j=1}^2 a_{ij}(x) \frac{\partial^2u^\e}{\partial x_i \partial x_j}(t,x) + b(x) \cdot \nabla u^\e(t,x), \hs4 u^\e(0,x) = g(x),
\end{equation}
where $g(x)$ is a bounded continuous function. 
\begin{theorem}\label{th_finite_pde_attractor}
Let $V_{12}$ and $V_{21}$ be defined by $\eqref{eq_finite_quasi}$ and assume $\ds \lim_{\e \to 0} T(\e) = \infty$. Let $g_1 := g(K_1)$ and
\begin{equation*}
	g_2 := \int_{K_2} g(x) \Big(|b(x)| \int_{K_2} \frac{dy}{|b(y)|} \Big)^{-1}dx.
\end{equation*}
If $x$ is attracted to $K_1$, then
\begin{equation*}
	\lim_{\e \to 0}u^\e(T(\e),x) = \begin{cases}
		g_1, & \text{if } \ds \limsup_{\e \to 0} \e \log T(\e) < V_{12} \text{ or } V_{12} > V_{21},
		\\ g_2, & \text{if } \ds \liminf_{\e \to 0} \e \log T(\e) > V_{12} \text{ and } V_{12} < V_{21}.
	\end{cases}
\end{equation*}
If $x$ is attracted to $K_2$, then
\begin{equation*}
	 \lim_{\e \to 0} u^\e(T(\e),x) = \begin{cases}
	 	g_2, & \text{if }\ds \limsup_{\e \to 0} \e \log T(\e) < V_{21} \text{ or } V_{21} > V_{12},
	 	\\ g_1, & \text{if } \ds \liminf_{\e \to 0} \e \log T(\e) > V_{21} \text{ and } V_{21} < V_{12}.
	 \end{cases}
\end{equation*}
\end{theorem}
The proof of Theorem \ref{th_finite_pde_attractor} uses Theorem \ref{th_finite_cycles} as well as the representation of the solution of problem $\eqref{eq_finite_Cauchy}$ in the form $u^\e(t,x) = \E_x g(X_t^\e)$, where the index $x$ in the expectation sign indicates that $X_0^\e = x$.

In the problems considered above, the set $\Me$ consists of a finite number of measures $\mu_1,...,\mu_{\tilde{\ell}}$. If the initial point $x = X_0^\e$ is attracted to one of these measures and the system is generic, then in the exponential time scale the limiting evolution consists of jumps from one point of $\Me$ to another one. Of course, for a fixed $\e > 0$ the evolution of measures is continuous, but the time spent outside small neighborhoods of the compacts $K_i$, is small compared to the time spent inside the neighborhoods. The existence of rough symmetries lead to a limiting motion in the convex envelope $\M$ of $\Me$. However, it turns out that for a certain class of perturbations, the limiting motion takes place in $\mathcal{M}$ even in the generic case. We demonstrate this phenomenon for the non-perturbed system shown in Figure \ref{fig_finite_example}. 

Consider the Cauchy problem
\begin{equation}\label{eq_finite_Cauchy_nonlinear}
	\frac{\partial u^\e}{\partial t}(t,x) =  L^\e u^\e(t,x), \hs4 u^\e(0,x) = g(x),
\end{equation}
where $L^\e$ is the nonlinear operator 
\begin{equation*}
L^\e v := \frac{\e}{2} \sum_{i,j=1}^2 a_{ij}(x,v) \frac{\partial^2v}{\partial x_i \partial x_j} + b(x) \cdot \nabla v.
\end{equation*}
The perturbation is now nonlinear. As before, assume $a_{ij}$ are twice continuously differentiable and uniformly positive definite. Moreover, for any $z \in \R$ we set
\begin{equation*}
	S_{0T}^z(\varphi) := \begin{cases}
		\ds \frac{1}{2}  \int_0^T \Big(a^{-1}(\varphi_s,z)(\dot{\varphi}_s - b(\varphi_s)) \Big)\cdot \Big(\dot{\varphi}_s - b(\varphi_s)\Big) ds, & \text{if }\varphi \in AC_{0T},
		\\ +\infty, & \text{otherwise},
	\end{cases}
\end{equation*}
and 
\begin{equation*}\label{eq_finite_quasi}
	V_{ij}(z) := \inf \{ S_{0T}^z(\varphi) : \varphi \in C_{0T}, \ \varphi_0 \in K_i, \ \varphi_T \in K_j, \ T \geq 0 \}.
\end{equation*}

\begin{figure}[h]
\centering
\includegraphics[scale = 1]{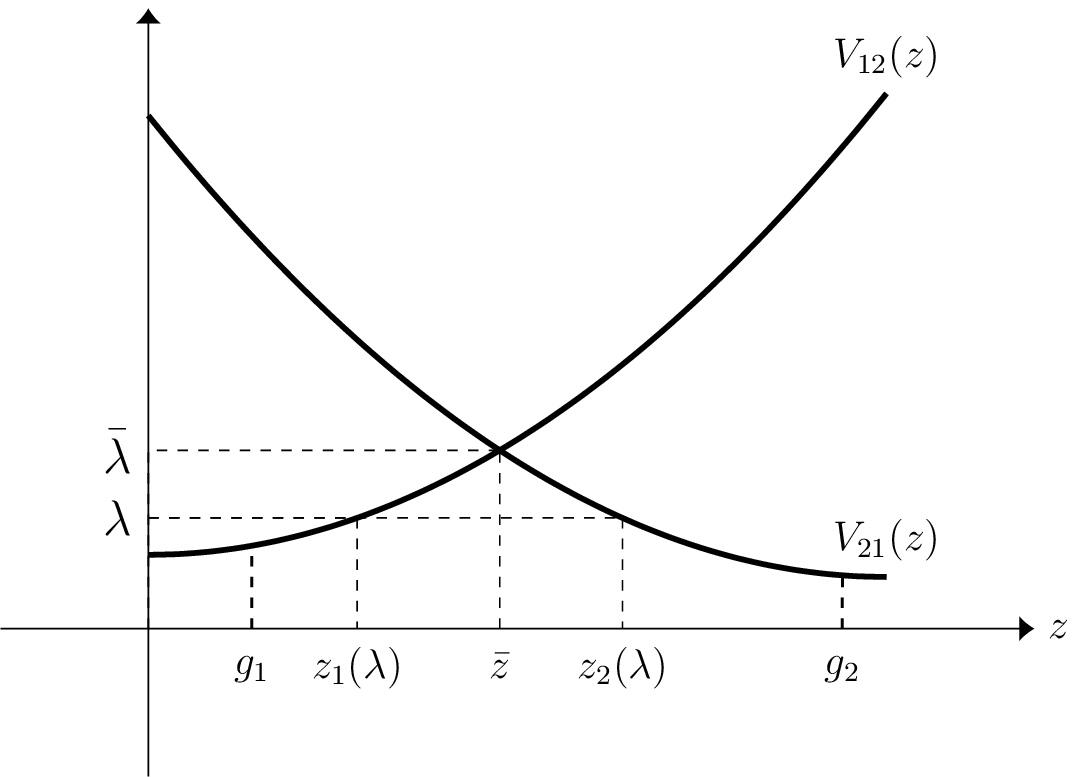}
\vspace*{-10mm}
\caption{}
\label{fig_finite_monotone}
\end{figure}

Suppose now that $V_{12}(z)$ and $V_{21}(z)$ are monotone as shown in Figure \ref{fig_finite_monotone}. Let $\bar{\lambda}$ and $\bar{z}$ be defined by the equality $V_{12}(\bar{z}) = V_{21}(\bar{z}) = \bar{\lambda}$. Additionally, let $z_1(\lambda)$ be the inverse function for $V_{12}(z) $ and $z_2(\lambda)$ be the inverse for $V_{21}(z)$, i.e. $V_{12}(z_1(\lambda)) = \lambda,$ $V_{21}(z_2(\lambda)) = \lambda)$. Assume $g_1$ and $g_2$ are as defined in Theorem \ref{th_finite_pde_attractor}. Without loss of generality, assume that $g_1 < g_2$.

One can consider the diffusion process associated with problem $\eqref{eq_finite_Cauchy_nonlinear}$. If $u^\e(t,x)$ is the solution of $\eqref{eq_finite_Cauchy_nonlinear}$, then define the random process $X_t^\e$ as the solution to the equation
\begin{equation}\label{eq_finite_Feynman_nonlin}
	 X_t^\e -x = \int_0^t \sigma(X_s^\e,u^\e(t-s,X_s^\e))dW_s + \int_0^t b(X_s^\e)ds,
\end{equation}
where $\sigma(x,u)$ is a Lipschitz continuous matrix-valued function such that $\sigma(x,u) \sigma^*(x,u)$ $ = a(x,u)$. The process $X_t^\e$ and function $u$ are then related by 
\begin{equation}\label{eq_finite_relate_Feynman}
	u^\e(t,x) = \E_x g(X_t^\e).
\end{equation}
Moreover, if one treats both $X_t^\e$ and $u^\e$ as unknowns, then the system defined by equations $\eqref{eq_finite_Feynman_nonlin}$ and $\eqref{eq_finite_relate_Feynman}$ is equivalent to the system $\eqref{eq_finite_Cauchy_nonlinear}$. This problem can be used to study problem $\eqref{eq_finite_Cauchy_nonlinear}$, as done in \cite{Freidlin85}.
\begin{theorem}
	Let $V_{12}(z)$ and $V_{21}(z)$ be defined above and as shown in Figure \ref{fig_finite_monotone}. Assume that $\ds \lim_{\e \to 0} T(\e) = +\infty$ and $X_0^\e = x \in \R^2$. Let $\mu_i$ be the $\d$-measure concentrated at $g_i$, $i \in \{1,2\}$. Then the following hold.
	
	If $\ds \liminf_{\e \to 0} \e \log T(\e) > \bar{\lambda}$, then the measures $\mu_{X_{T(\e)}^\e}$ converge weakly as $\e \downarrow 0 $ to the measure
	\begin{equation*}
	\frac{\bar{z} - g_1}{g_2 - g_1} \mu_2 + \frac{g_2 - \bar{z}}{g_2 - g_1} \mu_1.
	\end{equation*}
	If $x$ is attracted to $K_1$ and $\ds \lim_{\e \to 0} \e \log T(\e) = \lambda < \bar{\lambda}$, then the measures $\mu_{X_{T(\e)}^\e}$ converge weakly as $\e \downarrow 0 $ to the measure
	\begin{equation*}
		\frac{z_1(\lambda)-g_1}{g_2-g_1} \mu_2 + \frac{g_2 - z_1(\lambda)}{g_2-g_1} \mu_1.
	\end{equation*}
	If $x$ is attracted to $K_2$ and $\ds \lim_{\e \to 0} \e \log T(\e) = \lambda < \bar{\lambda}$, then the measures $\mu_{X_{T(\e)}^\e}$ converge weakly as $\e \downarrow 0 $ to the measure
	\begin{equation*}
		\frac{g_2 - z_2(\lambda)}{g_2-g_1} \mu_1 + \frac{z_2(\lambda)-g_1}{g_2-g_1} \mu_2.
	\end{equation*}
	Correspondingly ,
	\begin{equation*}
	\lim_{\e \to 0} u^\e(T(\e),x) = \begin{cases}
		\bar{z}, & \text{if } \ds \liminf_{\e \to 0} \e \log T(\e) > \bar{\lambda},
		\\ z_1(\lambda), & \text{if } \ds \lim_{\e \to 0 } \e \log T(\e) = \lambda < \bar{\lambda} \text{ and }x_0^\e = x \text{ is attracted to } K_1,
		\\ z_2(\lambda), & \text{if } \ds \lim_{\e \to 0 } \e \log T(\e) = \lambda < \bar{\lambda} \text{ and }x_0^\e = x \text{ is attracted to } K_2.
	\end{cases}
	\end{equation*}
\end{theorem}
To prove this theorem, note that $u^\e(t,x)$ is approximately constant for large $t$ and small $\e > 0$ inside the domain of attraction of each of the asymptotically stable compacts $K_i$. Therefore, the transition probabilities between neighborhoods of these compacts can be estimated using the action functional $S_{0T}^z$ with different $z$ in different domains of attraction. This allows us to apply arguments used in the linear case (\cite{Freidlin12}, Ch. 6). Problems similar to $\eqref{eq_finite_Cauchy_nonlinear}$ are considered in detail in \cite{Freidlin10}.

Until now we have considered non-degenerate white-noise-type perturbations of system $\eqref{eq_finite_dyn_sys}$. One can consider other types of perturbations that will lead to different limiting motions on $\Me$. Degenerate white-noise perturbations of $\eqref{eq_finite_dyn_sys}$ are of interest, in particular because of the Langevin equation
\begin{equation}\label{eq_finite_Langevin}
	\begin{cases}
	 \mu \ddot{q}(t) = b(q_t) - \lambda \dot{q}(t)  + \sqrt{\e} \dot{W}_t
	\\  q(0) = q \in R^n, \hs4 \dot{q}(0) = p \in \R^n.
	\end{cases}
\end{equation}
Here, $\mu$ and $\lambda$ are positive parameters, $W_t$ is the Wiener process in $\mathbb{R}^n$ and $b$ is a vector field in $\R^n$. If we rewrite $\eqref{eq_finite_Langevin}$ as a first order system in $\R^{2n}$, the noise will be degenerate. This leads to a different action functional,
\begin{equation*}
	S_{0T}^*(\varphi) = \frac{1}{2} \int_0^T \Big| \mu \ddot{\varphi}_s + \lambda \dot{\varphi}_s - b(\varphi_s) \Big|^2 ds,
\end{equation*}
which correspondingly leads to a different motion on $\Me$.

Instead of white-noise-type perturbations of $\eqref{eq_finite_dyn_sys}$, one can consider perturbations leading to Markov processes with jumps. The corresponding action functional is described in \cite{Freidlin12} (see also \cite{Wentzell91}).

One can also consider fast oscillating perturbations of $\eqref{eq_finite_dyn_sys}$, such as in the system
\begin{equation}\label{eq_finite_homogen}
	\dot{X}_t^\e = \tilde{b}(X_t^\e,h_{t/\e}), \hs4 X_0^\e = x \in \R^n,
\end{equation}
where $h_t$ is a stationary process with good enough mixing properties and is such that $\E \tilde{b}(x,h_t) = b(x).$ Then $X_t^\e$ converges to the solution of $\eqref{eq_finite_dyn_sys}$ on finite time intervals. However, on long time intervals growing together with $\e^{-1}$, $X_t^\e$ can have transitions between asymptotically stable attractors of the non-perturbed system. For instance, one can take as $h_t$ a non-degenerate diffusion process on a compact manifold. The action functional for the family of processes $X^\e_\cdot(t)$ defined by $\eqref{eq_finite_homogen}$ is calculated in \cite{Freidlin12}, Ch. 7. 

Let us now consider another example of a system satisfying Assumption \ref{assume1} with a finite number of ergodic probability measures. Suppose we have domains $D_1,...,D_m \subset \R^n$ and that each $D_k$ is homeomorphic to a ball and has a smooth boundary $\partial D_k$. Moreover, assume that $\partial D_i \cap \partial D_j = \emptyset$ for $i \neq j$. Denote by $\tilde {X}_t$ the diffusion process on $\R^n$ governed by the operator $L$,
\begin{equation*}
	L u(x) = \frac{1}{2} \sum_{i,j = 1}^n a_{ij}(x) \frac{\partial^2 u}{\partial x_i \partial x_j}(x) + \sum_{i=1}^n b_i(x) \frac{\partial u } {\partial x_i}(x).
\end{equation*}
We assume the coefficients of $L$ are bounded together with their first derivatives and the matrix $a(x) = (a_{ij}(x))$ is uniformly positive definite in $\R^n$. Moreover, let the process $\tilde{X}_t$ be positive recurrent (see, for instance, \cite{Freidlin85}). We now define $X_t$ to be the Markov process which coincides with $\tilde{X}_t$ outside of $\cup_{k=1}^m \partial D_k$, but when $X_t$ hits $\partial D_k$ it enters $D_k$ and stays there forever. We enforce this by imposing reflection boundary conditions on $\partial D_k$ in the direction of the interior co-normal $\eta^{-}(x)$. The precise definition of $X_t$ is given by its generator $(A,D(A))$ . The domain $D(A)$ contains the functions $u(x)$ that are continuous on $\R^n$, twice differentiable outside of $\cup_{k=1}^m \partial D_k$ and satisfy $\frac{\partial u(x)}{\partial \eta^{-}_k} \Big|_{\partial D_k} = 0$ for each $ k \in \{1,...,m\}$, while $Au(x) = Lu(x)$ outside of $\cup_{k=1}^m \partial D_k$ and $Lu(x)$ is continuous in $\R^n$. The process $X_t$ is our non-perturbed system. Ergodic invariant measures of $X_t$ are invariant measures of $X_t$ inside each of $D_k$ such that $D_k \cap \partial D_i = \emptyset$ for $i \neq k$. It is clear that Assumption \ref{assume1} is satisfied. In this case the perturbation consists of a small change in the gluing conditions on $\cup_{k=1}^m \partial D_k$. Namely, the perturbed gluing condition for the generator has the form
\begin{equation*}
	\Big( \e \frac{\partial u(x)}{\partial \eta^+_k} - \frac{\partial u(x)}{\partial \eta^-_k} \Big) \Big|_{\partial D_k} = 0,
\end{equation*} 
where $\frac{\partial}{\partial \eta^+}$ and $\frac{\partial}{\partial \eta^-}$ mean to calculate the derivative in the co-normal direction from outside and the inside of $D_k$ respectively. The perturbed process $X_t^\e$ associated to this generator can now exit the domains $D_k$ with small nonzero probabilities.

In this case, unlike in the previous problems, the measure $\mu_x$ corresponding to the initial point $x$ is not necessarily ergodic. If $\ds \lim_{\e \to 0}T(\e) = \infty$ but $\ds \lim_{\e \to 0} \e T(\e) = 0,$, then the distribution $X_{T(\e)}^\e$ will be close to $\mu_x$ as $\e \to 0$. However, if $\ds \lim_{\e \to 0} \e T(\e) > 0$, then $\mu_{X_{T(\e)}^\e}$ can make jumps from one point of $\M$ to another. In some cases, there will be a hierarchy of cycles. In general, there will be a hierarchy of Markov chains \cite{Freidlin17}, \cite{Freidlin20}.

Consider now a continuous time Markov chain $X_t$ with a finite phase space $\mathcal{E} = \{1,...,N\}$. As is well known, under mild assumptions such a chain has a finite set of invariant ergodic probability measures $\Me$. Moreover, for each $i \in \mathcal{E}$ the distribution of $X_t$ (with $X_0 = i$) converges as $t \to \infty$ to a probability measure $\mu_i$ belonging to the convex envelope $\M$ of $\Me$. Small perturbations of the chain can lead to transitions between different $\mu \in \M$, and in an appropriate time scale this motion on $\M$ can have a limit as $\e \downarrow 0$.

Let, for brevity, the non-perturbed Markov chain $X_t$ be the trivial chain that stays at the initial point forever with no transitions between states. In this case $\Me$ consists of $N$ unit measures: one at each state $i \in \mathcal{E}$. Let the perturbed chain $X_t^\e$ have the same phase space and a transition intensity matrix $Q^\e = \{q_{ij}^\e\}$ such that $\sum_{j=1}^N q_{ij}^\e = 0$ and $\lim_{\e \to 0} q_{ij}^\e = 0$. Put,
\begin{equation*}
	\mathbb{P}(X_{t+\d}^\e = j| X_t^\e = i\} = q_{ij}^\e \d +o(\d),
\end{equation*}  
as $\d \to 0$ for $i \neq j$. Assume also that $q_{ij}^\e > 0$ for each $i \neq j$ and $\e > 0$. We are interested in the time evolution of $\mu_{X_{T(t,\e)}^\e}$ on $\M$ and in the limit of $\mu_{X_{T(t,\e)}^\e}$ as $\e \to 0$ for an appropriate timescale $T(t,\e)$. Perturbations of Markov chains were considered by many authors in recent years (see for instance \cite{LandimArxiv} and the references therein, \cite{Freidlin14}, \cite{Freidlin17}).

\begin{figure}[h]
\centering
\includegraphics[scale = .9]{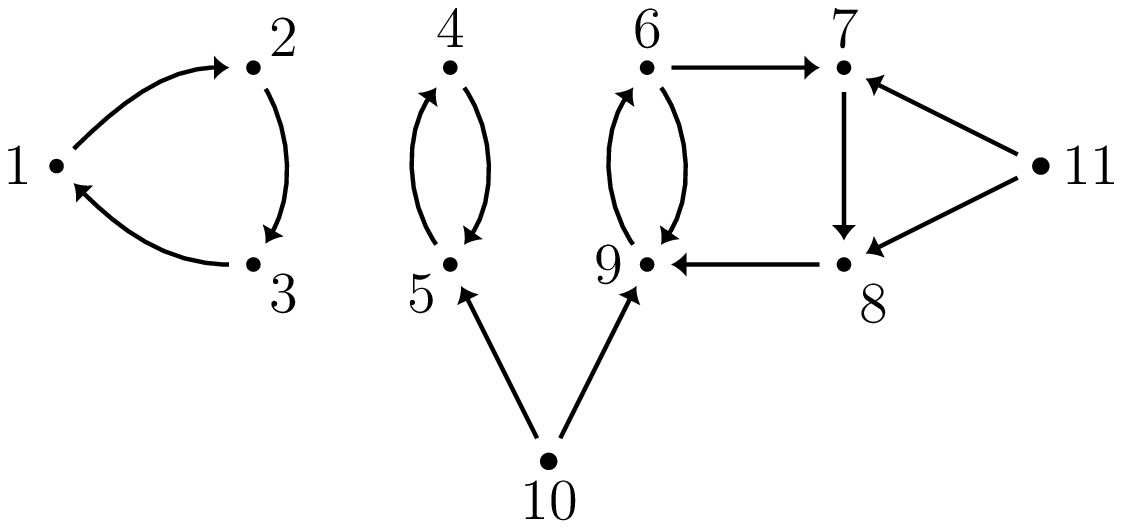}
\vspace*{-4mm}
\caption{}
\label{fig_finite_chain}
\end{figure}

Of course, without any additional assumptions one cannot expect regularity in the behavior of $\mu_{X_{T(t,\e)}^\e}$. Following \cite{Freidlin17}, we say that the perturbed chain is asymptotically regular if 
\begin{equation*}
\lim_{\e \to 0} \frac{q_{ij}^\e}{q_{mn}^\e} = \alpha_{ij;mn} \in [0,\infty]
\end{equation*}
exists for  all $i,j,m,n \in \mathcal{E}$ such that $i \neq j$ and $ m \neq n$. If the chain is asymptotically regular, one should then connect different points $i,j \in \mathcal{E}$ by an arrow if $\alpha_{i,j;i,k} > 0$ for all $k \in \mathcal{E} \setminus \{i\}$. It is then possible to uniquely choose disjoint sets $\mathcal{E}^{(1)},...,\mathcal{E}^{(m)} \subset \mathcal{E}$ such that the following conditions are satisfied.
\begin{enumerate}
	\item For any $k = 1,...,m$ and any two points $i,j \in \mathcal{E}^{(k)}$, there exists a sequence of arrows leading from $i$ to $j$.
	\item There are no arrows leading from each $\mathcal{E}^{(k)}$ to a point of $\mathcal{E} \setminus \mathcal{E}^{(k)}$.
	\item For each $k \in \mathcal{E} \setminus \Big( \bigcup_{i=1}^m \mathcal{E}^{(i)}\Big)$, there exists a sequence of arrows leading from $k$ to at least one of $\mathcal{E}^{(i)}$.
\end{enumerate}
An example of this decomposition is shown in Figure \ref{fig_finite_chain}. In the example, the sets are $\mathcal{E}^{(1)} = \{1,2,3\}$, $\mathcal{E}^{(2)} = \{4,5\}$ and $\mathcal{E}^{(3)} = \{6,7,8,9\}$, while points $10$ and $11$ are outside of $\cup_{i = 1}^3 \mathcal{E}^{(i)}$. Point $10$ is connected with $\mathcal{E}^{(2)}$ and $\mathcal{E}^{(3)}$, while point $11$ is connected with $\mathcal{E}^{(3)}$.

On each $\mathcal{E}^{(k)}$, consider a Markov chain $\mathcal{E}^{(k)}_\e$ with transition intensities 
\begin{equation*}
	q^\e_{ij} = \begin{cases}
		q_{ij}^\e & \text{if $i,j \in \mathcal{E}^k$ and $i\neq j$},
		\\   - \sum_{\ell \in \mathcal{E}^{(k)}, \ell \neq i} q_{i\ell}^\e & \text{if $i = j \in \mathcal{E}^k$}.
	\end{cases}
\end{equation*}
These Markov chains $\mathcal{E}^{(k)}_\e$ are called chains of rank $1$ or $1$-chains. One invariant ergodic distribution is defined on each chain $\mathcal{E}^{(k)}_\e$. Denote it by $\nu_k^\e$. One can show that there exists $\ds \lim_{\e \to 0}\nu^\e_k = \nu_k$.

By the time the process $X_t^\e$ does eventually leave $\mathcal{E}^{(k)}$, its distribution will be very close to $\nu_k^\e$. This allows us to calculate the asymptotics of transition intensities from $\mathcal{E}^{(k)}$ to $\mathcal{E}^{(i)}$ in an explicit form using $i$-graphs \cite{Freidlin12}. In this way we obtain a Markov chain of rank $2$. If the asymptotic regularity condition is satisfied for $2$-chains, one can introduce a chain of rank $3$ and so on until all $N$ states will be involved. One can give a simple complete asymptotic regularity condition which provides the construction of chains of all ranks. This construction is very similar to construction of the hierarchy of cycles mentioned above (see \cite{Freidlin12}, Ch.6). A cycle is, of course, a special type of a Markov chain. Having this hierarchy of chains, one can calculate the metastable distribution for each initial point and time scale in a similar fashion as in the case of cycles, as well as consider such phenomena as stochastic resonance, etc. \cite{Freidlin00}, \cite{Freidlin12}.

\begin{figure}[h]
\centering
\includegraphics[scale = .7]{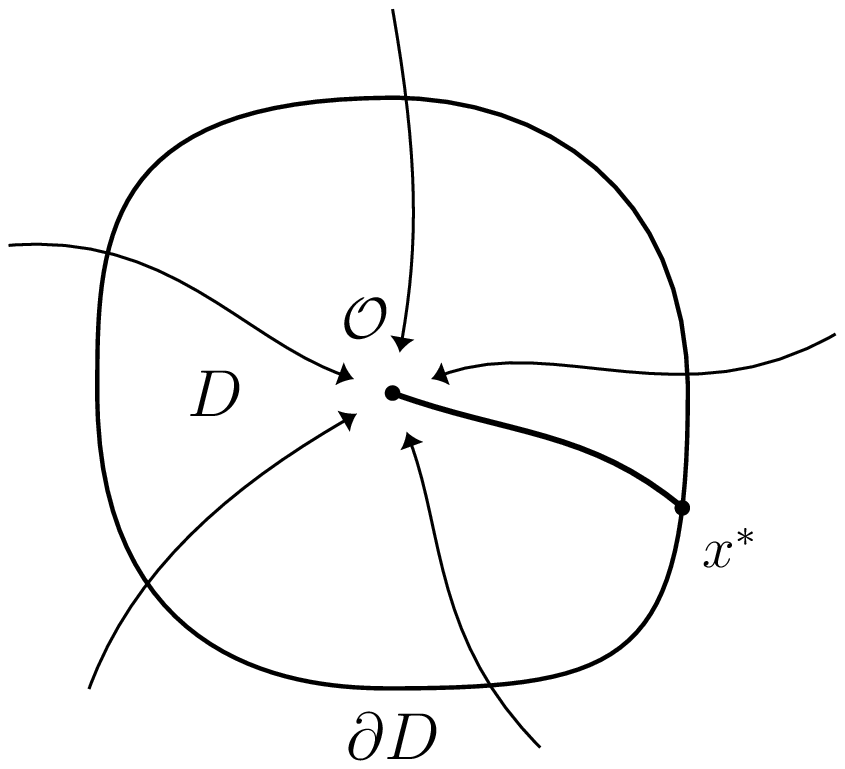}
\vspace*{-6mm}
\caption{}
\label{fig_finite_exit}
\end{figure}

Finally, let us mention the exit problem. Suppose we are given a bounded domain $D \subset \R^n$ with smooth boundary $\partial D$ and a vector field $b(x)$ satisfying for some $\alpha > 0$
\begin{equation*}
	 b(x) \cdot n(x) |_{x \in \partial D} \leq - \alpha.
\end{equation*}
Here $n(x)$ is the exterior normal to $\partial D$. This implies that trajectories of system $\eqref{eq_finite_dyn_sys}$ never leave the domain $D$ after entering it. For example, there could be an asymptotically stable equilibrium $\O$ inside $D$ that is attracting to the entire space $\R^n$, as in Figure \ref{fig_finite_exit}.

If small stochastic perturbations are added, the perturbed system $X_t^\e$ with $X_0^\e = x \in D$ may be able to leave $D$. For example, this is possible if $X_t^\e$ is the solution to $\eqref{eq_finite_stoc_pert}$ with a non-degenerate $a(x) = \sigma(x)\sigma^*(x)$ and bounded $|\beta(x)|$. Define the first exit time
\begin{equation*}
\tau^\e := \min \{t : X_t^\e \in \partial D\}.
\end{equation*}
What can one say about the asymptotics of $\tau^\e$, $\E_x \tau^\e$ and $X_{\tau^\e}^\e$ as $\e \to 0$? These questions concern the long-time behavior of the perturbed system; however, the non-perturbed system shown in Figure \ref{fig_finite_exit} has just one invariant probability measure: the $\d$-measure at point $\mathcal{O}$. In order to incorporate this into our approach, instead of $X_t^\e$ we consider the stopped process
\begin{equation*}
	\bar{X}_t^\e := \begin{cases}
		X_t^\e, \text{ if } t < \tau^\e,
		\\ X_{\tau^\e}^\e, \text{ if }t \geq \tau^\e.
	\end{cases}
\end{equation*}
It is clear that the first exit time $\tau^\e$ and the position at time $\tau^\e$ for $X_t^\e$ and $\bar{X}_t^\e$ are the same. The non-perturbed system $\bar{X}_t$ corresponding to $\bar{X}_t^\e$ is then defined as follows: $\bar{X}_t = X_t$ if $\bar{X}_0 = x \in D$, and $\bar{X}_t \equiv x $ if $\bar{X}_0 = x \in \partial D$. The system $\bar{X}_t$ has an invariant probability measure concentrated not just at $0$ but also at each point of the boundary $\partial D$ of the domain $D$. Therefore, the set $\Me$ for $\bar{X}_t$ can be parametrized by the points of the set $\{\mathcal{O}\} \cup \partial D$. One can prove (\cite{Freidlin12}, Ch.4, 6) that in the general case $\ds \lim_{\e \to 0} \e \log \tau^\e = V(x^*)$, where $V^* = \ds \min_{x \in \partial D} V(x)$ and $V(x)$ is the quasi-potential 
\begin{align*}
	V(x) & := \inf \Big\{ \ds \frac{1}{2}  \int_0^T \Big(a^{-1}(\varphi_s,z)(\dot{\varphi}_s - b(\varphi_s)) \Big)\cdot \Big(\dot{\varphi}_s - b(\varphi_s)\Big) ds: \\ & \varphi \in C_{0T}, \  \varphi_0 = 0, \ \varphi_T = x, \ T > 0 \Big\},
\end{align*}
and $\ds \lim_{\e \to  0} X_{\tau^\e}^\e = x^*$ if $x^*$ is the unique minimum of $V(x)$ on $\partial D$. 
\section{Systems with $\Me$ homeomorphic to a graph}\label{sec_graph}
Consider a Hamiltonian system with one degree-of-freedom
\begin{equation}\label{eq_graph_Hamil_unper}
	\dot{X_t} = \bar{\nabla} H(X_t), \hs4 X_0 = x = (x_1,x_2) \in \R^2.
\end{equation}
We assume that the Hamiltonian $H(x)$ is smooth, $\ds \lim_{|x| \to \infty} H(x) = \infty$ and $\ds \min_{x \in \R^2} H(x) = 0$. Moreover, suppose $H(x)$ has a finite number of critical points and all of them are non-degenerate. Assume for simplicity that $H(x)$ has a different values at each critical point.
\begin{figure}[h]
\centering
\includegraphics[scale = .9]{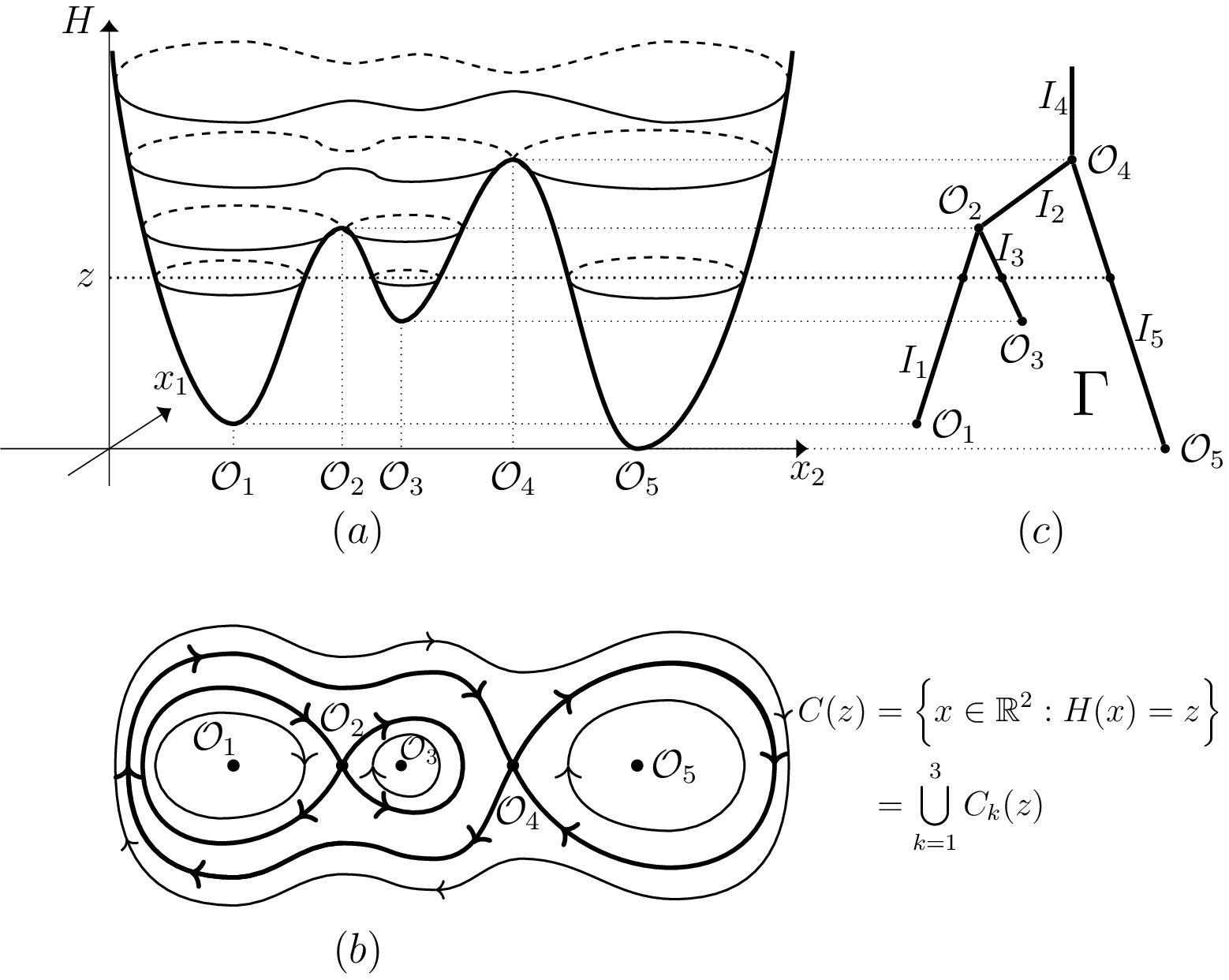}
\vspace*{-10mm}
\caption{}
\label{fig_graph_main}
\end{figure}

Let $\Gamma$ be the graph which ``counts'' the connected components of all level sets of $H(x)$ (Figure \ref{fig_graph_main}). Set
\begin{equation*}
	C(z) := \{x \in \R^2:H(x) = z \} = \bigcup_{k} C_k(z),
\end{equation*} 
where $C_k(z)$ are the individual connected components of $C(z)$, whose number depends on $z$. Exterior vertices of $\Gamma$ ($\mathcal{O}_1$, $\mathcal{O}_3$, $\mathcal{O}_5$ in Figure \ref{fig_graph_main}) correspond to the minima of $H(x)$. Each interior vertex ($\mathcal{O}_2$, $\mathcal{O}_4$ in Figure \ref{fig_graph_main}) corresponds to an eight-shaped curve ($\gamma_1,\gamma_2$ in Figure \ref{fig_graph_main}). If one numbers the edges of $\Gamma$ ($I_1,...I_5$ in Figure \ref{fig_graph_main}), each point $y \in \Gamma$ can be characterized by its value of $H$ and by the number of the edge containing it. In this way one can introduce on $\Gamma$ a global coordinate system: $(h,i)$ is the point on the edge number $i$ such that $H(x) = h$ on the corresponding level set component (for details see Ch. 8 of \cite{Freidlin12}). 

Let $Y:\R^2 \to \Gamma$ denote the mapping $ x \mapsto  (H(x),i(x))$ and consider the following perturbed Hamiltonian system:
\begin{equation}\label{eq_graph_Hamil_pert}
	\dot{\tilde{X}}_t^\e = \bar{\nabla}H(\tilde{X}_t^\e) + \e \beta( \tilde{X}_t^\e) + \sqrt{\e} \sigma(\tilde{X}_t^\e) \dot{W}_t, \hs4 \tilde{X}_0^\e = x.
\end{equation}
The generator of the process $\tilde{X}_t^\e$ is $\tilde{L}^\e = \e L^\e$, where
\begin{equation*}
L^\e u(x) = \frac{1}{\e} \bar{\nabla} H(x) \cdot \nabla u + \beta(x) \cdot \nabla u + \frac{1}{2} \sum_{i,j=1}^2 a_{ij}(x) \frac{\partial^2 u }{\partial x_i \partial x_j}.
\end{equation*}
As usual, $a(x) = \sigma(x) \sigma^*(x)$ is assumed to be uniformly positive definite. Note that the operator $L^\e$ is the generator of the process $X_t^\e := \tilde{X}_{t/\e}^\e$. 

First, suppose the Hamiltonian $H(x)$ has just one well. In this case the graph $\Gamma$ consists of just one edge $I=\{h \geq 0\}$. Hence, one can omit the integer-valued coordinate so that the mapping $Y$ is just given by $Y(x) = H(x)$. The classical averaging principle in this case says (for instance, see Ch. 8 of \cite{Freidlin12}) that $Y_t^\e := Y(X_t^\e) = H(X_t^\e) $ converges weakly to the diffusion process $Z_t$ on $I$ governed by the operator $\bar{L}$ defined by
\begin{equation}\label{eq_graph_Lbar}
	\bar{L}u(z) = \frac{1}{2} \bar{a}(z) \frac{d^2u}{dz^2} + \bar{\beta}(z) \frac{du}{dz}.
\end{equation}
Here $\bar{a}(z) = \frac{1}{T(z)} \int_{G(z)} \mathrm{div} \Big( a(x) \nabla H(x) \Big) dx$ and $\bar{\beta}(z) = \frac{1}{T(z)} \int_{G(z)} \mathrm{div} \beta(x)dx,$ where $T(z) = \oint_{C(z)} \frac{d\ell}{|\nabla H(x)|}$ is the period of rotation along $C(z) = \{x:H(z) = z\}$ of the non-perturbed system and $G(z)$ is the domain in $\R^2$ bounded by $C(z)$. The vertex $\mathcal{O}\in I$ is inaccessible for the process $Z_t$.

If the Hamiltonian $H(x)$ has more than one well, as in Figure \ref{fig_graph_main}, the graph $\Gamma$ will have interior vertices. Inside each edge $I_i \subset \Gamma$, up to the time of first exit from $I_i$, one can describe the limit of $Y_t^\e$ as $\e \to 0$ in a similar way as for the case of one well. However, in this case it turns out that the interior vertices can be accessible for $Z_t$ in a finite time. Therefore, the classical averaging principle should be supplemented with a description of what the trajectory of the limiting process on $\Gamma$ should do after hitting an interior vertex.

One can expect that the limiting process $Z_t$ on $\Gamma$ is a continuous Markov process. All continuous Markov processes on graphs, with some mild regularity conditions, were described in \cite{Freidlin12}.

Next, let $(A,D_A)$ be the generator $A$ of the process $Z_t$ together with its domain $D_A$. Let $u:\Gamma \to \R$ be a continuous on $\Gamma$ and smooth inside the edges function that also belong to $D_A$. Then inside each edge $I_k \in \Gamma$, $Au(z)$ coincides with $\bar{L}_k u(z)$, defined as in $\eqref{eq_graph_Lbar}$, and at each interior vertex $\mathcal{O}_i$ the following gluing condition should be satisfied:
\begin{equation}\label{eq_graph_glue}
	\alpha_i A u(\mathcal{O}_i) = \sum_{j : I_j \sim \mathcal{O}_i} \pm \gamma_{ij} D_j u(\mathcal{O}_i).
\end{equation}
Here $\alpha_i, \gamma_{ij}$ are non-negative constants and $D_ju(\mathcal{O}_i)$ is the derivative of $u(z)$ (where $z = (h,j)$) in $h$ calculated along the edge $I_j$, while the notation $I_j \sim \mathcal{O}_i$ means that $I_j$ is attached to $\mathcal{O}_i$. Note that $Au(\mathcal{O})$ is the common value of $\bar{L}_j u(\mathcal{O}_i)$ for each $j$ such that $I_j \sim \mathcal{O}_i$. Signs $+$ ($-$) should be taken if the coordinate $h$ on $I_j$ is greater (smaller) than $H(\mathcal{O}_i)$. The operators $\bar{L}_k$ together with the gluing conditions define the process $Z_t$ on $\Gamma$ in a unique way.

The operators $\bar{L}_k$ can be calculated by standard averaging using $\eqref{eq_graph_Lbar}$. How can one find the constants $\alpha_i, \gamma_{ij}$ in $\eqref{eq_graph_glue}$? To answer this, first let the perturbations in $\eqref{eq_graph_Hamil_pert}$ be such that $L^\e = \frac{1}{\e} \bar{\nabla}H \cdot \nabla +R$ where $R$ is any formally self-adjoint second order operator. It is easy to check that the Lebesgue measure $\Lambda$ in $\R^2$ is invariant for $X_t^\e$ for each $\e \geq 0$. This implies that the projection $\lambda$ of $\Lambda$ onto $\Gamma$ defined by the mapping $Y$ is invariant for $Y_t^\e$ and for their limit as $\e \to 0$. It turns out that there exists just one set of constants $\alpha_i, \gamma_{ij}$ for which the limiting process on $\Gamma$ has the prescribed invariant measure $\lambda$; since $\lambda(\{\mathcal{O}_i\}) = 0$, $\alpha_i = 0$. 

Now, if in addition to a self-adjoint 2nd order term, $L^\e$ also has a term with first derivatives, then the measure in the space of trajectories on a finite time interval corresponding to this operator is absolutely continuous with respect to the measure without this new term. This implies that addition of the term with first derivatives will not change the gluing conditions. Since each second order elliptic operator with smooth coefficients can be written as a sum of the self-adjoint part and first order operator, one can find the constants $\alpha_i,\gamma_{ij}$ in the general case. The complete proof is available in Ch. 8 of \cite{Freidlin12}.

One can consider other types of perturbations of $\eqref{eq_graph_Hamil_unper}$. The simplex $\M$ will still be the same, as well as graph $\Gamma$, but the limit of $Y(X_t^\e)$ as $\e \to 0$ may be different. Moreover, this limit may not exist without an appropriate regularization. Consider the classical case of pure deterministic perturbations: $\sigma(x) \equiv 0$ in $\eqref{eq_graph_Hamil_pert}$. Suppose for brevity that $\mathrm{div} \beta(z) < 0$. If $H(x)$ has just one well then $\Gamma$ consists of one edge $I = \{h\geq 0 \},$ and $Y_t^\e = Y(\tilde{X}_{t/\e}^\e) = H(\tilde{X}_{t/\e}^\e)$ converges uniformly on each finite time interval as $\e \downarrow 0$ to the solution of the averaged equation 
\begin{equation*}
	\dot{Z}_t = \bar{\beta}(Z_t), \hs4 Z_0 = H(\tilde{X}_0^\e).
\end{equation*}
Next suppose that $H(x)$ has more than one well so that $H(x)$ has at least one saddle point $\O_2$ (as in Figure \ref{fig_graph_main} for instance). Then the trajectory $Z_t$ starting from $z_0 \in (H(\mathcal{O}_2),H(\mathcal{O}_4))$ hits $H(\mathcal{O}_2)$ in a finite time $t(z_0)$. After $t(z_0)$, $\tilde{X}_{t/\e}^\e$ will go to one of two wells separated by the saddle point. These wells will alternate as $\e \to 0$ so that $\ds \lim_{\e \to 0} Y_t^\e$ for $t > t(z_0)$ will not exist. 

There are many ways to regularize this problem. For instance, one can add a small random perturbation to the initial point: instead of a fixed initial point $x$ one can take a random point distributed uniformly on the $\d$- neighborhood of $x$. Then the solution of $\eqref{eq_graph_Hamil_pert}$ with $\sigma(x) \equiv 0$ will be a stochastic process $\tilde{X}_t^{\e,\d}$. One can consider first the limit of $Y_t^{\e,\d} = Y(\tilde{X}_{t/\e}^{\e,\d})$ as $\e \to 0$, and then subsequently the $\d \to 0$ limit in the sense of weak topology. If the Hamiltonian has just one saddle point, such a double limit $Y_t$ exists for any perturbation $\beta(x)$ that is sufficiently regular. The process $Y_t$ will be a deterministic motion inside each edge and at the saddle points it will have some stochasticity. For instance, in the case of the saddle point $\mathcal{O}_2$ in Figure \ref{fig_graph_main}, $Y_t$ will go to the edges $I_1$ and $I_3$ with probabilities proportional to $\int_{G_1} | \mathrm{div} \beta(x) | dx$ and $\int_{G_2} |\mathrm{div} \beta(x) | dx$, respectively. Here $G_1$ and $G_2$ are the left and right domains bounded by the curve $\gamma_1$ in Figure 
\ref{fig_graph_main}. If $H(x)$ has more than one saddle point the double limit $\ds \lim_{\d \to 0} \lim_{\e \to 0} Y(\tilde{X}_{t/\e}^{\e,\d})$ does not exist for a generic class of perturbations. Nonetheless, one can always regularize the problem by addition a small white-noise type perturbation to the equation. Let now $\tilde{X}_t^{\e,\d}$ be the solution of $\eqref{eq_graph_Hamil_pert}$ with $\sigma(x)$ replaced by $\sqrt{\d} \sigma(x)$ for some $ 0 < \d \ll 1.$ Then the double limit $\ds \lim_{\d \to 0} \lim_{\e \to 0} Y( \tilde{X}_{t/\e}^{\e,\d})$ exists for any perturbation $\beta(x)$. The limiting process will be the same for various matrices $\sigma(x)$ with non-degenerate $a(x) = \sigma(x) \sigma^*(x)$, and in addition this limit will coincide with the one obtained by perturbation of the initial point, if the latter exists. This means that the stochasticity displayed in the limit of the deterministically perturbed system is actually an intrinsic property of the deterministic system, which are close to one-degree-of-freedom Hamiltonian systems if the latter have saddle points (see \cite{Freidlin12}, Ch. 8). 

Equation $\eqref{eq_graph_Hamil_pert}$ with a degenerate matrix $a(x) = \sigma(x) \sigma^*(x)$ arises naturally in certain physical models including, for instance, when random perturbations of an oscillator are considered. This case is studied in \cite{Freidlin98}.

Consider now a flow in $\R^3$ defined by the generalized Landau-Lifshitz equation
\begin{equation}\label{eq_graph_landau}
	\dot{X}_t = \nabla F(X_t) \times \nabla G(X_t), \hs4 X_0 = x \in \R^3.
\end{equation}
Here $F(x)$ and $G(x)$, $x \in \R^3$, are smooth generic functions with $F(x) \geq 0$ and $\ds \lim_{|x| \to \infty} F(x) = \infty$. The flow $X_t$ is incompressible since
\begin{equation*}
	\mathrm{div} (\nabla F \times \nabla G) = \nabla G \cdot \mathrm{curl} \nabla F - \nabla F \cdot \mathrm{curl} \nabla G = 0.
\end{equation*} 
Both $F(x)$ and $G(x)$ are integrals of motion, since for example,
\begin{equation*}
	\frac{dF(X_t)}{dt} = \nabla F(X_t) \cdot (\nabla F(X_t) \times \nabla G(X_t)) = 0.
\end{equation*} 
The classical Landau-Lipschitz equation has $F(x) = \frac{|x|^2}{2}$. This equation arises in magnetization theory. From the physical perspective, perturbations that preserve the first integral $F(x)$ are of interest. Thus one should consider the flow $X_t$ and its perturbations on the $2D$-surface $S_z = \{x \in \R^3:F(x) = z\}$. It is convenient in this case to write the stochastic part of the perturbation in the
Stratonovich form. Assume that $S_z$ consists of one connected component and $|\nabla F(x)| \neq 0$ for each $x \in S_z$. As is known, under these assumptions the flow $X_t$ on $S_z$ has an invariant measure with the density $\frac{c}{|\nabla F(x)|}$ with respect to the area element on $S_z$ for an appropriate constant $c$. In particular, if $F(x) = \frac{|x|^2}{2}$ then this measure is the uniform distribution on the sphere $S_z$.

In the case of equation $\eqref{eq_graph_landau}$, one can consider the graph $\Gamma$ counting the connected components of the level sets of $G(x)$ on $S_z$. One ergodic probability measure is concentrated on each connected component so that $\Gamma$ is homeomorphic to the set $\Me$ provided with the topology of weak convergence. Each compact orientable 2D-manifold topologically can be characterized by its genus $\chi$, a non-negative integer. If $\chi = 0 $, the manifold is a topological sphere. If the genus of $S_z$ is $0$ and white-noise-type or pure deterministic perturbations preserving $S_z$ are considered, then the calculation of the limiting process on $\Me$ (or equivalently on $\Gamma$) is very similar to the case of perturbations of a one-degree-of-freedom Hamiltonian system \cite{Freidlin11}. One should just keep in mind that the measure with the density $\frac{c}{|\nabla F(x)|}$ on $S_z$ is invariant when calculating the constants in gluing conditions $\eqref{eq_graph_glue}$ so that all $\alpha_k$ will be equal to $0$. 

The situation is different if the genus of $S_z$ is greater than $0$. The structure of trajectories of area-preserving systems on compact 2D-manifolds was described in \cite{Dolgopyat08}. These results imply that the set of invariant ergodic probability measures for such systems can be parametrized by a graph $\Gamma$, but some of the vertices, which we will call heavy, correspond to measures that have a density with respect to the area on $S_z$. The number of heavy vertices is equal to the genus $\chi$. In our case, $X_t$ does not in general preserve the area on $S_z$ but the existence of an invariant measure with a continuous positive density implies that the structure of invariant measures will be the same. The arguments mentioned above allow us to calculate the limiting motion inside the edges and the constants in the gluing conditions $\eqref{eq_graph_glue}$. The coefficients $\alpha$ at heavy vertices will be positive. In order to complete the proof, one should prove that the limit process in this case also will be Markovian. The complete proof was given in \cite{Dolgopyat08}.

If pure deterministic perturbations of $\eqref{eq_graph_landau}$ that preserve $S_z$ are considered, one can regularize the system. The limiting motion in this case also will also be deterministic inside the edges with a stochastic branching at interior vertices. But, in addition it will spend a random exponentially distributed time at the heavy vertices. The parameters of these exponential times are also independent of the regularization \cite{Dolgopyat12}.

Consider now the diffusion process $X_t$ in $\R^n$ governed by the operator $L_0$,
\begin{equation*}
L_0 u(x) = \frac{1}{2} \mathrm{div} (a(x) \nabla u) + b(x) \cdot \nabla u, \hs4 x \in \R^n.
\end{equation*}
We assume that the matrix $a(x) = (a_{ij}(x)) $ is non-negative definite with bounded first and second derivatives entries, and $b(x)$ is a Lipschitz continuous vector field. Then one can show there exists a Lipschitz continuous matrix $\sigma(x)$ and vector field $\tilde{b}(x)$ such that $X_t$ is the solution of
\begin{equation*}
\dot{X}_t = \sigma(X_t) \dot{W}_t + \tilde{b}(X_t), \hs4 X_0 = x \in \R^n,
\end{equation*}
where $W_t$ is the Wiener process in $\R^n$. We say that a smooth function $H(x)$ is a first integral for the process $X_t$, if $\mathbb{P}_x(H(X_t) = H(x)) = 1$ for each $x \in \R^n$.

It is easy to check that $H(x)$ is a first integral for $X_t$ if $a(x) \nabla H(x) \equiv 0$ and $b(x) \cdot \nabla H(x) \equiv 0$. Assume these hold and also that $H(x)$ has a finite number of critical points which are all non-degenerate. Moreover, assume $\min_{x \in \R^n} H(x) = 0$ and $\lim_{|x| \to \infty} H(x) = \infty$. Assume that $\lambda_1(x) |e|^2 \leq (a(x) e \cdot e) \leq \lambda_2(x) |e|^2$ if $e \perp \nabla H(x)$ for $x \in \R^n$. Here $\lambda_1(x) > 0$ if $\nabla H(x) \neq 0$, $\lambda_2(x) < K < \infty$ inside a ball containing all critical points of $H(x)$; if $\nabla H(\O) = 0$ then $\lambda_1(x) \geq k_1 |x - \O|^2$ and $\lambda_2(x) \leq k_2|x-\O|^2$ for $x$ in a neighborhood of $\O$. Then the process $X_t$ has one ergodic probability measure $M_{z,i}$ on each connected component $C_i(z)$ of $C(z) = \{x \in \R^n:H(x) = z\}$. If $C_i(z)$ contains no critical points, $M_{z,i}$ has a density $m_{z,i}(x)$ on $C_i(z)$. Otherwise, $M_{z,i}$ is the $\d$-measure at the critical point belonging to $C_i(z)$. The set $\Me$ in this case can e parametrized by the graph $\Gamma$. As before, one can number the edges of $\Gamma$, introduce the coordinates $(z,k)$ and the projection $Y:\R^n \to \Gamma, Y(x) = (z(x),k(x))$. This process $X_t$ in $\R^n$ is our non-perturbed system. 

To describe the perturbed process, define the operator $L_1$, 
\begin{equation*}
	L_1u = \frac{1}{2} \mathrm{div} \Big( a_1(x) \nabla u(x) \Big) + \beta(x) \cdot \nabla u(x), \hs4 x \in \R^n.
\end{equation*}
Assume the entries of the matrix $a_1(x)$ and the vector field $\beta(x)$ are bounded and smooth enough and that $L_1$ is a strictly elliptic operator. Denote by $\tilde{X}_t^\e$ the process governed by $\tilde{L}^\e = L_0 + \e L_1$; $\tilde{X}_t^\e$ is the perturbed process. Next set $X_t^\e = \tilde{X}_{t/\e}^\e$ and $Y_t^\e = Y(X_t^\e)$. 

If the critical points of $H(x)$ are non-degenerate, then in a neighborhood of each critical point $\mathcal{O}_k$ one can write $H(x) = x_1^2 + ...+x_r^2 - x_{r+1}^2 - ... x_n^2$ in an appropriate coordinate system. If $\mathcal{O}_k$ is an extremum of $H(x)$, then $r = 0$ or $r = n$, while if $\O_k$ is a saddle point then $1 \leq r \leq n - 1$. Moreover, if $\mathcal{O}_k$ is an extremum of $H(x)$, just one edge is attached to the vertex corresponding to $\mathcal{O}_k$ in the graph $\Gamma$. One can prove that if $r > 1$ and $n - r > 1$ (this can happen if $n > 3$) exactly two edges are attached to the corresponding vertex. If $r = 1$ or $n - r = 1$, then a small neighborhood of $\mathcal{O}_k$ will be divided into 3 parts, but two of these parts can be united outside of the small neighborhood. Therefore, if $r = 1$ or $r = n-1$, two or three edges are attached to the vertex corresponding to $\mathcal{O}_k$. This can happen if $n \geq 4$.

Similar to the case of one-degree-of-freedom Hamiltonian systems, the evolution of $\mu_{X_t^\e}$ can be described by the process $Y_t^\e = Y(X_t^\e)$ on $\Gamma$. It is proved in \cite{Freidlin04} that $Y_t^\e$ converges weakly on each finite time interval as $\e \to 0$ to the diffusion process $Z_t$ on $\Gamma$, which is defined by its generator $A$. The domain $D_A$ of $A$ contains continuous functions $u:\Gamma \to \R$ that are bounded and smooth inside the edges such that $u'(z)$ is continuous at all interior vertices $\mathcal{O}_k$ with just two edges are attached; if three edges are attached to $\mathcal{O}_k$, a gluing condition of form $\eqref{eq_graph_glue}$ with $\alpha_k = 0$ should be satisfied at $\mathcal{O}_k$. Inside each edge $I_j$, $Au(z) = \bar{L}_j u(z)$, where 
\begin{equation*}
	\bar{L}_j u(z) = \frac{1}{2} \bar{a}_j(z) u''(z) + \bar{\beta}_j(z) u'(z).
\end{equation*}
Here the coefficients $\bar{a}_j(z)$ and $\bar{\beta}_j(z)$ are calculated by averaging. Explicit formulas for $\bar{a}_j$ and $\bar{\beta}_j$ as well as for the gluing conditions are available in \cite{Freidlin04}. Note that many important characteristics of the limiting process $Z_t$ on $\Gamma$ can be calculated explicitly as solutions of simple ordinary differential equations.

Consider now the Dirichlet problem for the operator $\tilde{L}^\e = L_0 + \e L_1$ in a bounded domain $G$ with smooth boundary $\partial G$. Here $L_0$ and $L_1$ are operators described above. For the sake of brevity, we assume that the vector field $b(x)$ defines an incompressible flow so that the Lebesgue measure is invariant for the diffusion process $X_t$ governed by $L_0$. Let $u^\e$ be the solution to the Dirichlet problem
\begin{equation}\label{eq_graph_dirich}
\begin{cases}
\tilde{L}^\e u^\e(x)= 0, & x \in G,
\\ 	 u^\e(x)  = \psi(x), & x \in \partial G,
\end{cases}
\end{equation}
where $\psi(x)$ is a continuous function on $\partial G$.

As discussed in Section \ref{sec_finite}, in order to apply our approach to the Dirichlet problem, one should consider the stopped process
\begin{equation*}
	\bar{X}_t^\e = \begin{cases}
		X_t^\e, & \text{if } t < \tau^\e,
		\\ X_{\tau^\e}^\e, & \text{if } t \geq \tau^\e,
	\end{cases}
\end{equation*}
where $\tau^\e = \inf \{t > 0:X_t^\e \in \partial G \}$ and $X_t^\e$ is the process in $\R^n$ governed by $\tilde{L}^\e$. The corresponding stopped non-perturbed process is then given by
\begin{equation*}
\bar{X}_t = \begin{cases}
	X_t^0, & \text{if } t \leq \tau^0,
	\\ X_{\tau^0}^0, & \text{if }t \geq \tau^0.
\end{cases}
\end{equation*}
Denote by $\pi(x,\gamma)$ the distribution on $\partial G$ of $X_{\tau^0}^0$ given $X_0^0 = x$. The process $\bar{X}_t^0$ satisfies Assumption \ref{assume1}, but the measure $\mu_x$ attracting an initial point $x \in \R^n$ is not necessarily ergodic. For instance, it can be $\pi(x,\cdot)$.

Let $\Gamma$ be the graph corresponding to $H(x)$ and $Y: \R^2 \to \Gamma$ be the corresponding mapping. Assume that each level set contains no more than one critical point and that $Y(\partial G)$ does not contain any vertex of $\Gamma$. Denote by $\Theta$ the set of $z = (h,i) \in \Gamma$ such that $C_i(h) \cap \partial G \neq \emptyset$. Then if $Y(x) \in \Theta$, the process $ \bar{X}_t$ given $\bar{X}_0 =x$ approaches the measure $\pi(x,\cdot)$ as $t \to \infty$. The measure $\pi(x,\cdot)$ belongs to the simplex $\M$ of invariant probability measures, but it is not ergodic if the support of $\pi(x,\cdot)$ consists of more than one point. The set of ergodic probability measures $\Me$ of $\bar{X}_t$ can be parametrized by $\partial G \cup (\Gamma \setminus \Theta)$.
\begin{figure}[h]
\centering
\includegraphics[scale = .8]{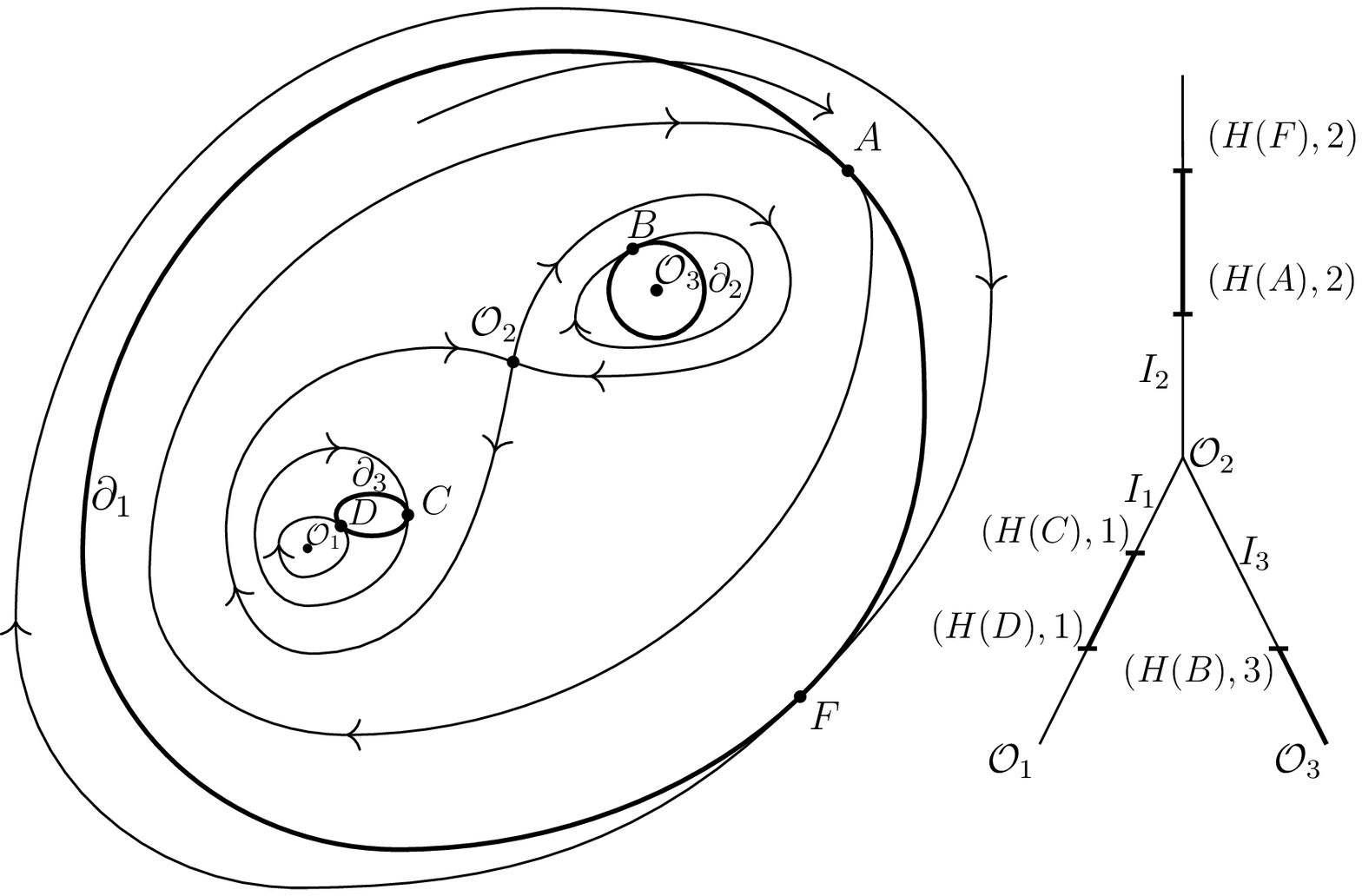}
\vspace*{-14mm}
\caption{}
\label{fig_graph_hamil}
\end{figure}

In the example shown in Figure \ref{fig_graph_hamil}, $\partial G$ consists of three parts $\partial_1$, $\partial_2$ and $\partial_3$; the set $\Theta$ consists of three closed intervals $[(H(A),2), (H(F),2)]$ on the edge $I_2$, $[(H(D),1), (H(C),1)]$ on $I_1$ and $[\mathcal{O}_3, (H(B),3)]$ on $I_3$.

Assume that the measure-valued function $\pi(x,\cdot)$ depends continuously on $x \in G \cup \partial G$ in the topology of weak convergence. Let the set $\Theta$ consist of a finite number of connected pieces. For example, in Figure \ref{fig_graph_hamil}, $\Theta$ is the union of the three intervals mentioned above. The boundary $\partial \Theta$ of $\Theta$ consists of exterior vertices and points $(h,k) \in \Gamma$ such that $C_k(h) \cap \partial G \neq \emptyset$. Assume that the set $C_k(h) \cap \partial G$ consists just of one point (points $A,B,C,D$ in Figure \ref{fig_graph_hamil}) if $(h,k) \in \partial \Theta$. Then $\psi(Y^{-1}(h,k))$ is correctly defined. 

Define a function $\bar{u}: G \cup \partial G \to \R$,
\begin{equation*}
	\bar{u}(x) = \begin{cases} \ds
		\int_{\partial G} \psi(y) \pi(x,dy), & \text{if }Y(x) = (h,k) \in \Theta \setminus \partial \Theta,
		\\  \psi(Y^{-1}(h,k)), & \text{if } (h,k) \in \partial \theta, 
		\\  v_k(h), & \text{if }(h,k) \notin \theta.
	\end{cases}
\end{equation*} Here, $v_k$ is the solution to 
\begin{equation*}
	\bar{L}_k v_k(h) = \frac{1}{2} \bar{a}_k(h) v_k''(h) - \bar{b}_k(h) v_k'(h) = 0, \hs4  v_k(h)|_{\partial \Theta} = \psi(Y^{-1}(h,k)),	
\end{equation*}
such that $v_k(h)$ as well as $\bar{L}_k v_k(h)$ are continuous on $\Gamma \setminus \Theta$ and at each interior vertex $\mathcal{O}_k$ such that just two edges are attached to $\O_k, v_k(h)$ and $v_k'(h)$ are continuous. If three edges are attached to $\O_k$, then the following gluing condition is satisfied:
\begin{equation*}
\sum_{i:I_i \sim \O_k} \pm \beta_i(\O_k) D_iv_i(h)|_{(h,i)=\O_k} = 0,
\end{equation*}
with sign $+$ if $h > H(\O_k)$ on $I_i$ and $-$ otherwise. Above, we are using
\begin{align*}
& \bar{a}_k(h) = \oint_{C_k(h)} \frac{a_1(x) \nabla H(x) \cdot \nabla H(x)}{|\nabla H(x)|}ds,
\\ & \bar{\beta}_k(h) = \oint_{C_k(h)} \frac{L_1H(x)}{|\nabla H(x)|} ds,
\\ & \beta_i(\O_k) = \oint_{\tilde{C}_i(\O_k)} \frac{a_1(x) \nabla H(x) \cdot \nabla H(x)}{|\nabla H(x)|}ds,
\end{align*}
with $\tilde{C}_i(\O_k) = \{ x \in \R^n:Y(x) = \O_k\} \bigcap \partial \{x \in \R^n: Y(x) \in I_i\}$. 

These conditions define $\bar{u}: G\cup \partial G \to \R$ in a unique way. Moreover $\bar{u}(x)$ can be calculated explicitly.

\begin{theorem}
	Let $u^\e(x)$ be the solution of problem $\eqref{eq_graph_dirich}$. If the above mentioned conditions are satisfied then for any $x \in D$,
	\begin{equation*}
		\lim_{\e \to 0} u^\e(x) = \bar{u}(x).
	\end{equation*}
\end{theorem}
The proof of this result is based on the representation of solution $u^\e(x)$ in the form $u^\e(x) = \E_x \psi(X_{\tau^\e}^\e)$ and on Theorem 8.2.1 from \cite{Freidlin12} which described the limiting evolution of $\mu_{X_{t/\e}^\e}$ as $\e \to 0$.
\begin{figure}[h]
\centering
\includegraphics[scale = .9]{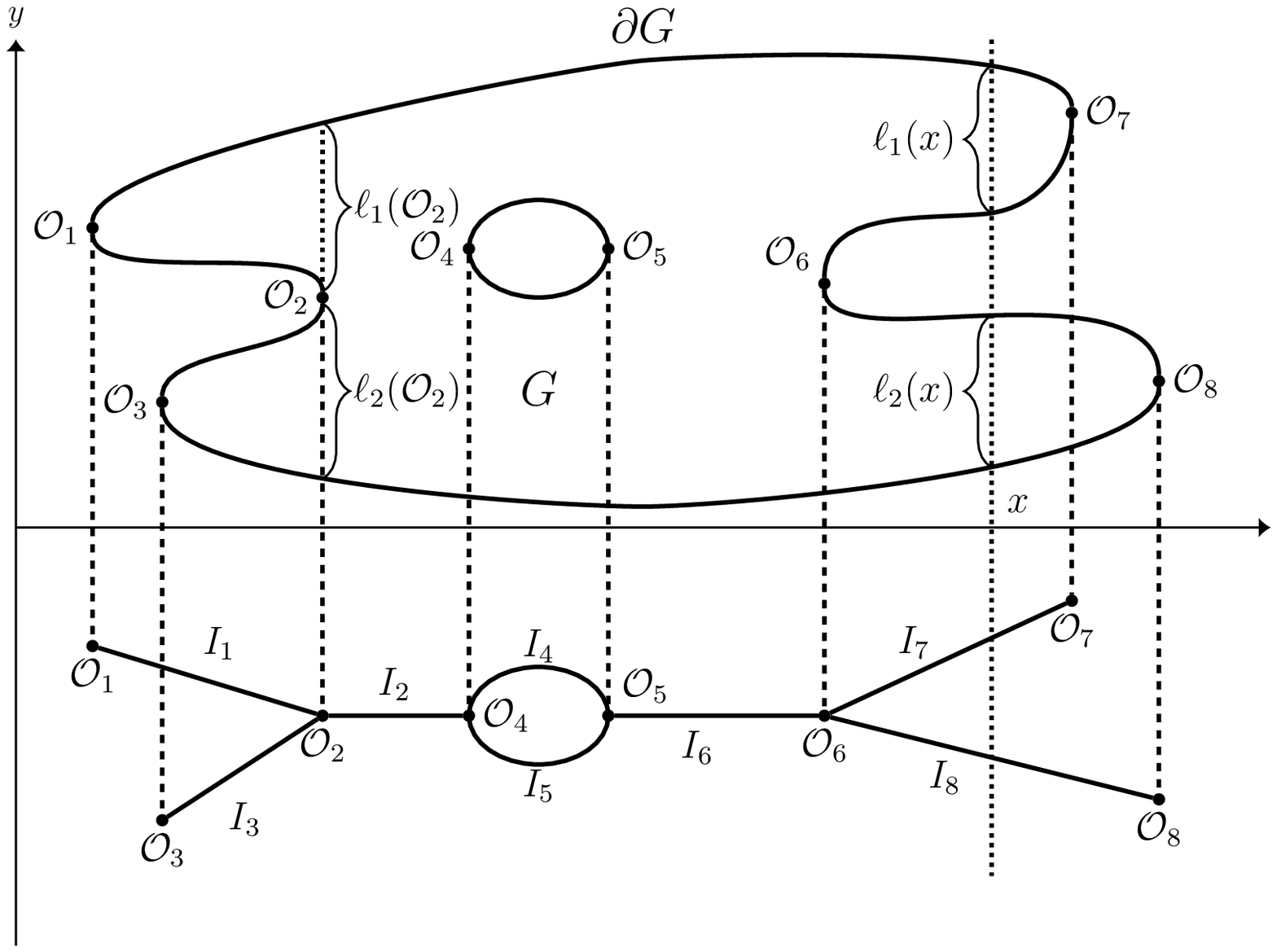}
\vspace*{-10mm}
\caption{}
\label{fig_graph_2Dreg}
\end{figure}

Finally, we mention the Neumann problem for a PDE with a small parameter.  a Let $G$ be a two-dimensional closed region such as that shown in Figure \ref{fig_graph_2Dreg}. Consider the Neumann problem
\begin{equation}\label{eq_graph_channel}
\begin{cases} \ds
	L^\e u^\e(x,y) := \frac{1}{2} \frac{\partial^2 u^\e}{\partial x^2} + \frac{1}{2\e} \frac{\partial^2 u^\e}{\partial y^2} = f(x,y), & (x,y) \in G
	\\ \ds  \frac{\partial u^\e(x,y)}{\partial n_\e(x,y)} = 0  & (x,y) \in \partial G .
	\end{cases}
\end{equation}
Here $n_\e(x,y)$ is the interior co-normal to $\partial G$ corresponding to the operator $L^\e$; i.e. $n_\e(x,y) = (\e n_1(x,y),n_2(x,y))$, where $(n_1(x,y),n_2(x,y))$ is the normal to $\partial G$. This type of problem arises when one considers the diffusion in narrow channels or layers. As is well known, problem $\eqref{eq_graph_channel}$ is solvable for $\e > 0$ if $\int_G f(x,y)dxdy = 0 $. The solution of $\eqref{eq_graph_channel}$ is unique up to an additive constant, so that to single out a unique solution one can assume that $\int_G u(x,y) = 0$.

Now consider the diffusion process $(X_t^\e,Y_t^\e)$ defined by
\begin{equation}\label{eq_graph_2Dproc}
\begin{cases}
	\ds \dot{X}_t^\e = \sqrt{\e} \dot{W}_t^1 + \e n_1(X_t^\e,Y_t^\e) \dot{\phi}_t^\e, & (X_0^\e, Y_0^\e) = (x,y) \in G,
	\\ \ds \dot{Y}_t^\e = \dot{W}_t^2 + n_2(X_t^\e,Y_t^\e) \dot{\phi}_t^\e, & \phi_0^\e = 0,
	\end{cases}
\end{equation}
where $(W_t^1,W_t^2)$ is the Wiener process in $R^2$ and $\phi_t^\e$ is the local time on $\partial G$ (see, for instance \cite{Freidlin85}). The solution of the Neumann problem can then be written as follows (see \cite{Freidlin85}):
\begin{equation*}
u^\e(x,y) = - \int_0^\infty \E_{x,y} f(X_{s/\e}^\e, Y_{s/\e}^\e)ds.
\end{equation*}
Therefore, to understand the limiting behavior of $u^\e(x,y)$ as $\e \to 0$, one should consider long-time evolution of the process defined by $\eqref{eq_graph_2Dproc}$.

The process $(X_t^\e,Y_t^\e)$ can be considered as a perturbation of the following two-dimensional process $(X_t,Y_t)$:
\begin{equation}\label{eq_graph_2dproc_unper}
\begin{cases}
	 \dot{X}_t = 0, & (X_0,Y_0) = (x,y),
	\\  \dot{Y}_t = \dot{W}_t^2 + n_2(X_t,Y_t) \dot{\phi}_t, & \phi_0 = 0.
	\end{cases}
\end{equation}
The set of ergodic invariant probability measures of the process $\eqref{eq_graph_2dproc_unper}$ in $G$ can be parametrized by the graph $\Gamma$ shown in Figure \ref{fig_graph_2Dreg}: the vertices of $\Gamma$ correspond to the points $(x,y) \in G$ where the tangent line to $\partial G$ is vertical. We assume that there are a finite number of such points. One or two invariant measures correspond to each point $x$ from the projection of $G$ on the $x$-axis. For instance, the intersection of the vertical line at the point $(x,0)$ in the picture with $G$ consists of two intervals, $\ell_1(x)$ and $\ell_2(x)$. Each of these intervals supports one ergodic probability measure of the non-perturbed system - the uniform distribution on this interval. The density of such a measure on an interval $\ell_i(x)$ is  $\frac{1}{\ell_i(x)}$ (we denote by $\ell_i(x)$ the interval itself and its length). Let us number the edges of $\Gamma$ and consider the mapping $Y:G \to \Gamma$. One can prove (see \cite{Freidlin12Neumann}) that $Y(X_{t/\e}^\e, Y_{t/\e}^\e)$ converges weakly as $\e \downarrow 0$ on each finite time interval $ [0,T]$ to a diffusion process $(\bar{X}_t,k_t)$ on $\Gamma$ which is defined by its generator $A$ and domain $D_A$ as follows:
\begin{enumerate}[(i)]
\item In the interior of an edge $I_i$, 
\begin{equation*}
Au(x,i) = \bar{L}_i u(x,i) := \frac{1}{2 \ell_i(x)} \frac{ d }{dx} \Big( \ell_i(x) \frac{du}{dx} \Big).
\end{equation*}
\item The domain $D_A$ contains functions $u:\Gamma \to \R$ that are continuous, smooth inside the edges such that $L_iu(x,i)$ is continuous on $\Gamma$, and satisfy the following conditions at the vertices. If three edges $I_{i_1},$ $ I_{i_2}$ and $ I_{i_3}$ are attached to $\O_k$ (points $\O_2,\O_4,\O_5,\O_6$ in Figure \ref{fig_graph_2Dreg}) then $\O_k$ separates two intervals $\ell_1(\O_k)$ and $\ell_2(\O_k)$ and the following gluing condition should be satisfied at $\O_k$:
\begin{equation*}
	\ell_1(\O_k) D_{i_1}u(\O_k) + \ell_2(\O_k) D_{i_2}u(\O_k) = (\ell_1(\O_k) + \ell_2(\O_k)) D_{i_3} u(\O_k),
\end{equation*}
where $D_i$ means differentiating along $I_i$. If $\O_k$ is an exterior vertex and $I_i \sim \O_k$ then $\lim_{x \to \O_k} \ell_i(x) D_i u(x) = 0$.
\end{enumerate}
Define a function $\bar{f}:\Gamma \to \R$,
\begin{equation*}
	\bar{f}(x,k) = \frac{1}{\ell_k(x)} \int_{\ell_k(x)} f(x,y)dy.
\end{equation*}
Then one can derive from \cite{Freidlin12Neumann} that the solution $u^\e(x,y)$ satisfying $\int_G u^\e(x,y)dxdy = 0$ converges as $\e \downarrow 0$ uniformly in $G$ to the function $v(Y(x,y)) = v(x,k(x,y))$ where $v:\Gamma \to \R$ satisfies the equation 
\begin{equation*}
Av(x,k) = \bar{f}(x,k), \hs4
\sum_{i : I_i \subset \Gamma} \int_{I_i} v(x,i) \ell_i(x)dx = 0.
\end{equation*}
The function $v(x,i)$ can be calculated explicitly.

\section{Open book as a phase space of the long-time evolution}\label{sec_book}
Roughly speaking, an open book space is a set consisting of a finite number of pieces of $n$-dimensional manifolds (pages) glued together at manifolds of dimension less than $n$ (the binding of the book). A graph is an example of an open book: edges glued at the vertices. Open book spaces arise naturally as the set homeomorphic to the collection of ergodic probability measures $\Me$ of a system. However, in many interesting examples $\Me$ itself may not be homeomorphic to an open book, but rather an essential part of $\Me$ will be homeomorphic to an open book space. This allows one to describe the long-time evolution of the system as a motion on the open book.

\begin{figure}[h]
\centering
\includegraphics[scale = .7]{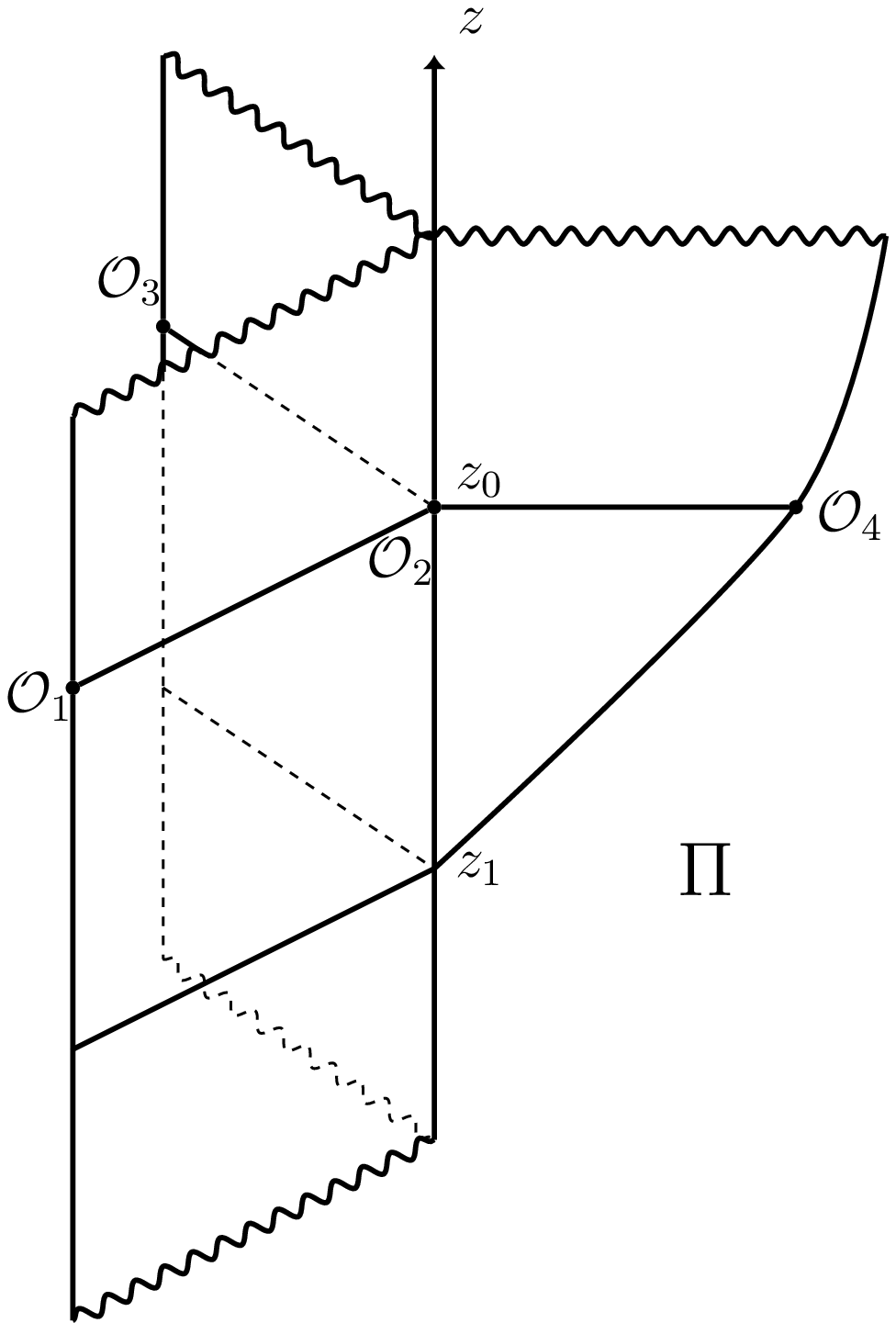}
\vspace*{-6mm}
\caption{}
\label{fig_book_lose}
\end{figure}

For example, in the case of Landau-Lipschitz equation $\eqref{eq_graph_landau}$, if the perturbations do not destroy the first integral $F(x)$, then the whole motion happens on $2D$-surface $S_z = \{x \in \R^3: F(x) = z\}$. The set $\Me^z$ of ergodic probability measures concentrated on $S_z$ can be parametrized by a graph $\Gamma_z$ ``counting'' the measures just on $S_z$. More precisely, the graph $\Gamma_z$ counts connected components of the level sets of $G(x)$ on $S_z$. If perturbations do not preserve $S_z$, then all ergodic probability measures of system $\eqref{eq_graph_landau}$ should be considered. It is clear that the set $\Me$ of all such measures is the union of $\Me^z$ for various $z$. For instance, if for $z = z_0$ the corresponding graph $\Gamma_{z_0}$ has three edges $\{\O_1\O_2, \O_3\O_2, \O_4\O_2\}$ as in Figure \ref{fig_book_lose}, when $z$ changes the graph can lose an edge (if $z < z_1$ in Figure \ref{fig_book_lose}). Conversely, new pages can appear as well. The whole set $\Me$ can thus be parametrized by the full open book $\Pi$. 

Consider another example where the set $\Me$ is homeomorphic to an open book. Define the operator $L^\e$, 
\begin{equation*}
	L^\e u(x) := \frac{1}{2} \sum_{i,j = 1}^n a_{ij}\Big(x,\frac{x}{\e}\Big) \frac{\partial^2 u}{\partial x_i \partial x_j} + \sum_{i=1}^n b_i \Big(x,\frac{x}{\e}\Big) \frac{\partial u}{\partial x_i},
\end{equation*}
where we assume first that the coefficients $a_{ij}(x,y)$ and and $b_i(x,y)$ are bounded, smooth functions from $\R^n\times\R^n \to \R$ and $1$-periodic in $y_1,...,y_n$. Assume the matrix $(a_{ij}(x,y))$ is uniformly positive definite. We are interested in the following homogenization problems: (1) what is the limiting behavior as $\e \downarrow 0$ of the diffusion process $X_t^\e$ in $\R^n$ governed by $L^\e$ and (2) what is the limiting behavior of the solutions of various boundary and initial-boundary problems related to the operator $L^\e$?

The diffusion process governed by $L^\e$ can be described by the equation 
\begin{equation*}
	\dot{X}_t^\e = b\Big(X_t^\e, \frac{X_t^\e}{\e}\Big) + \sigma \Big( X_t^\e, \frac{X_t^\e}{\e}\Big) \dot{W}_t, \hs4  X_0^\e = x \in \R^n,
\end{equation*}
where $b(x,y) = (b_z(x,y),...,b_n(x,y))$, $\sigma(x,y) \sigma^*(x,y) = a(x,y) = (a_{ij}(x,y))$ and $W_t$ is the $n$-dimensional Wiener process.

Next, put $\tilde{X}_t^\e=  X_{\e^2t}^\e$ and $\tilde{Y}_t^\e = \e^{-1} X_{\e^2t}^\e$. The pair $(\tilde{X}_t^\e, \tilde{Y}_t^\e)$ form (a degenerate) $2n$-dimensional process satisfying the equations
\begin{equation}\label{eq_book_degen_pair}
\begin{cases}
	 \dot{\tilde{X}}_t^\e = \e^2 b(\tilde{X}_t^\e,\tilde{Y}_t^\e) + \e \sigma(\tilde{X}_t^\e, \tilde{Y}_t^\e) \dot{W}_t,
	\\  \dot{\tilde{Y}}_t^\e = \e b(\tilde{X}_t^\e, \tilde{Y}_t^\e) + \sigma(\tilde{X}_t^\e, \tilde{Y}_t^\e)\dot{W}_t.
	\end{cases}
\end{equation}
System $\eqref{eq_book_degen_pair}$ is a perturbation of the following system:
\begin{equation}\label{eq_book_nonper}
\begin{cases}
	 \dot{\bar{X}}_t = 0, & \bar{X}_0 = x,
	\\ \ds  \dot{\bar{Y}}_t = \sigma(\bar{X}_t, \bar{Y}_t) \dot{W}_t, & \bar{Y}_0 = \ds \Big\{ \frac{x}{\e}\Big\}, 	
\end{cases}
\end{equation}
where $\{\frac{x}{\e}\}$ is the vector of fractional parts of $\frac{x}{\e}$ and the variable $\bar{Y}_t$ changes on the $n$-dimensional unit torus $\mathbb{T}^n$. For each $x \in \R^n$, the process $\bar{Y}_t^x$ defined by equation $\dot{\bar{Y}}_t^x = \sigma(x,\bar{Y}_t^x) \dot{W}_t$ on $\mathbb{T}^n$ has a unique invariant probability measure $\mu_x$ which is ergodic. The collection of such measures $\mu_x$, $x \in \R^n$ is the set $\Me$ for the non-perturbed system $\eqref{eq_book_nonper}$. In this case $\Me$ can be parametrized by points of $\R^n$. Note that each $\mu_x$ has a density $m_x$ in $\mathbb{T}^n$ which is the unique solution of the equation $\sum_{i,j=1}^n \frac{\partial^2}{\partial y_i \partial y_j} (a_{ij}(x,y) m_x(y)) = 0, y \in \mathbb{T}^n$, such that $\int_{\mathbb{T}^n} m_x(y)dy = 1$. Here $x \in \R^n$ is a parameter. 

Define the operator
\begin{equation*}
	\bar{L}u(x) = \frac{1}{2} \sum_{i,j = 1}^n \bar{a}_{ij}(x) \frac{\partial^2u}{\partial x_i \partial x_j} + \sum_{i=1}^n \bar{b}_i(x) \frac{\partial u}{\partial x_i},
\end{equation*}
where $\ds \bar{a}_{ij}(x) = \int_{\mathbb{T}^n} a_{ij}(x,y) m_x(y) dy$ and $\ds \bar{b}_i(x) = \int_{\mathbb{T}^n} b_i(x,y) m_x(y)dy.$ It was shown in \cite{Freidlin64} that the process $X_t^\e = \tilde{X}_{t/\e^2}^\e$ converges weakly as $\e \downarrow 0$ on each finite time interval to the process $\bar{X}_t$ in $\R^n$ (which is the parametrization of $\Me$) governed by the operator $\bar{L}$.

This result holds for a non-degenerate matrix $a_{ij}$. Actually, this result can be slightly improved. Instead of the non-degeneracy assumption, one can assume that the process defined by $\eqref{eq_book_nonper}$ has a unique invariant measure for each initial point $x \in \R^n$. However if there exists an open set $G \subset \R^n$ such that the process for $x \in G$ has more than one invariant probability measure, the homogenization can be non-complete and should be modified. We describe the situation in a special case. In particular, to avoid difficulties unrelated to the homogenization, we consider differential operator with discontinuous coefficients.

Let $a(y)$ be a smooth $1$-periodic function such that $a(z_1) = a(z_2) = a(z_3) = 0$ for $0 < z_1  < z_2 < z_3 < 1$ and $a(y) > 0$ if $y \notin \{z_1,z_2,z_3\}.$ Let $b(x_1,y)$ be a bounded smooth $1$-periodic in $y$ function. Put 
\begin{equation*}
	a\Big(x_1,\frac{x_2}{\e}\Big) := \begin{cases} 1, & \text{if } x_1 < 0,
	\\ a(x_2/\e), & \text{if } x_1 > 0, \end{cases}
\end{equation*}
and consider the operator $L^\e$ which for $\{(x_1,x_2) : x_1 \neq 0 \}$ is defined as follows:
\begin{equation*}
L^\e u(x_1,x_2) = \frac{1}{2} \frac{\partial^2u } {\partial x_1^2} + \frac{1}{2} a (x_1, x_2/\e) \frac{\partial^2 u}{\partial x_2^2} + b(x_1, x_2/\e) \frac{ \partial u} { \partial x_1}.
\end{equation*}
The operator $L^\e$ is considered on the domain of bounded continuous functions $u(x,y)$ that are smooth inside of each half-plane $\{x_1 < 0 \} $ and $\{x_1 > 0\}$ and such that $\frac{ \partial u}{\partial x_1}$ and $L^\e u$ are continuous in $\R^2$. Let $(X_1^\e(t),X_2^\e(t))$ be the diffusion process in $\R^2$ governed by $L^\e$ (together with its domain). One can check that such a process exists and is unique for each $\e > 0$.
\begin{figure}[h]
\centering
\includegraphics[scale = .7]{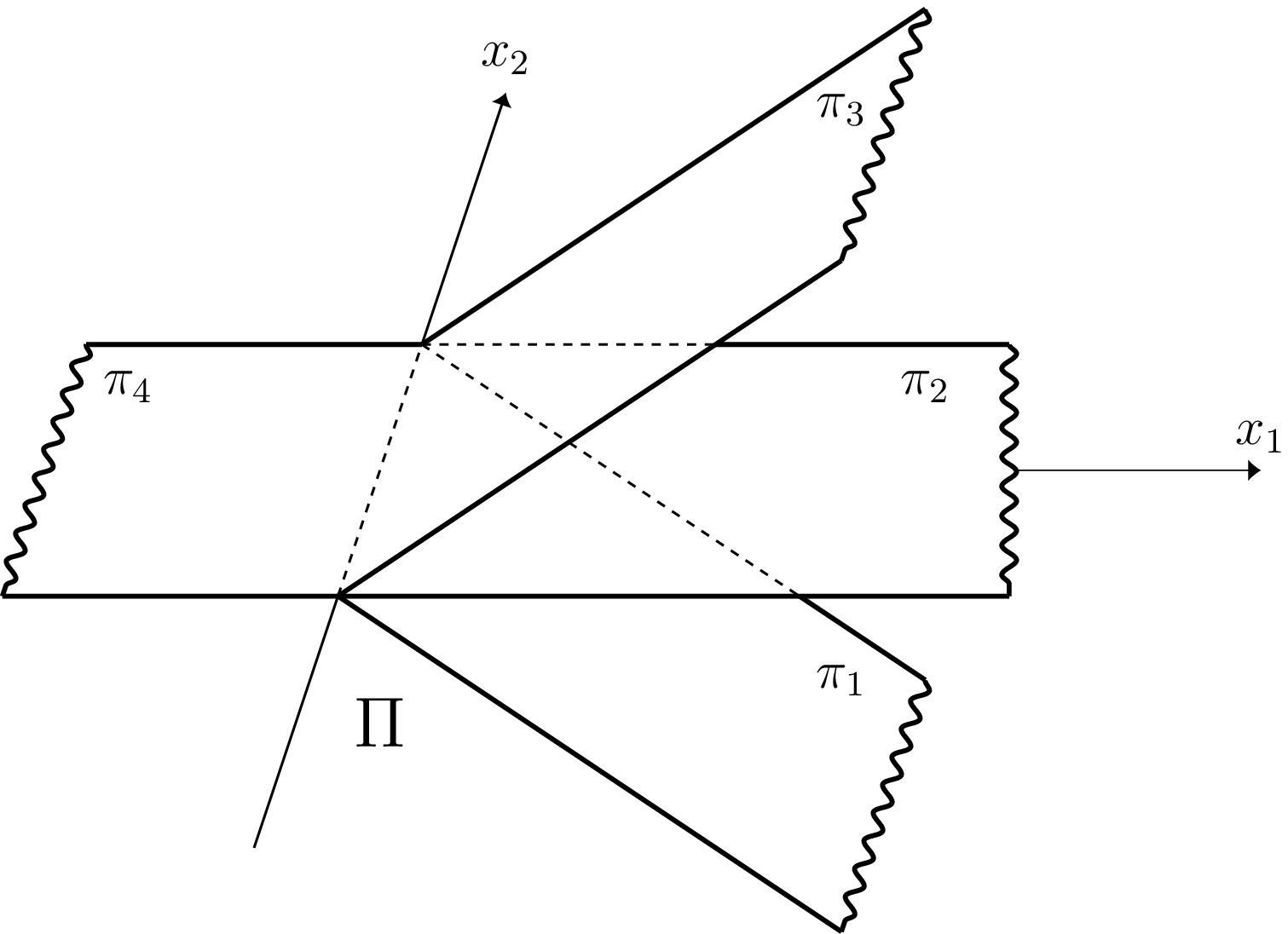}
\vspace*{-4mm}
\caption{}
\label{fig_book_homog}
\end{figure}

Let $\Pi$ be the open book space shown in Figure \ref{fig_book_homog}. It has $4$ pages $\pi_1$, $\pi_2$, $\pi_3$ and $\pi_4,$ and the binding. We provide $\Pi$ with the coordinate system $(x_1,x_2,k)$. Consider the diffusion process $Y_t$ on $\Pi$ defined as follows. On each page $\pi_i$, define an operator
\begin{equation*}
	\bar{L}_i = \frac{1}{2} \frac{\partial^2 u}{\partial x_1^2} + \frac{1}{2} \bar{a}_i \frac{\partial^2u}{\partial x_2^2}  + \bar{b}_i(x_1) \frac{\partial u}{\partial x_1},
\end{equation*}
where $\bar{a}_i = 0$ and $ \bar{b}_i(x_1) = b_i(x_1,z_i)$ for $i = 1,2,3$ while $\bar{a}_4 = 1$ and $\ds \bar{b}_4(x_1) = \int_0^1 b(x_1,y)dy$.

Now let $Y_t = (X_1(t), X_2(t),k_t)$ be the diffusion process on $\Pi$ governed by operator $\bar{A}$ defined as follows. Let $\bar{A}$ coincide with $\bar{L}_i$ inside each page $\pi_i$, and let its domain contain continuous on $\Pi$ and smooth inside the pages functions $u(x_1,x_2,k)$ such that $L_iu$ is continuous on $\Pi$ and  the following gluing condition is satisfied on the binding $\{x_1 = 0\}$:
\begin{align*}
	\Big( 2 D_4 u(x_1,x_2) & - (1+z_1-z_2) D_1u(x_1,x_2)  - (z_3-z_1) D_2 u(x_1,x_2) 
	\\ & - (1+z_2-z_3) D_3 u(x_1,x_2) \Big) \Big| _{x_1 = 0} = 0.
\end{align*}
Here $D_iu$ means differentiation in $x_1$ on the page $\pi_i$.

Note that the non-perturbed process related to $(X_1^\e(t),X_2^\e(t))$ is defined by equation $\eqref{eq_book_nonper}$. If $x_1 < 0$ then the non-perturbed trajectory starting at $(x_1,x_2)$ approaches the invariant measure $\mu_{x_1,x_2}$, which is the uniform distribution on $\{(x_1,x_2)\} \times \mathbb{T}$. If $x_1 > 0$ and $\{\frac{x_2}{\e}\}= \tau_\e$ is situated between neighboring $z_i,z_j \in \mathbb{T}$, then the invariant measure corresponding to $(x_1,x_2) \in \R^2$ is 
\begin{equation*}
\mu_{x_1,x_2} = \Big|\frac{z_i - \tau_\e}{z_i-z_j}\Big| \d_{x_1,x_2,z_j} + \Big|\frac{z_j - \tau_\e}{z_i-z_j}\Big| \d_{x_1,x_2,z_i},
\end{equation*} 
where $\d_{x_1,x_2,z}$ is the $\d$-measure concentrated at $(x_1,x_2,z) \in \R^2 \times \mathbb{T}$ where $\mathbb{T}$ is the unit circle. This implies that $(X_1^\e(t),X_2^\e(t))$ has no limit as $\e \downarrow 0$ and the system should be regularized. Note that for $x_1 > 0$ the measure $\mu_{x_1,x_2}$ is not in general ergodic.

We regularize the system by adding a random perturbation to the initial point. Denote by $(X_1^{\e,\d}(t), X_2^{\e,\d}(t))$ the process $(X_1^\e(t),X_2^\e(t))$ with the initial point $(x_1,x_2 + \d \xi)$ where $\d > 0$ and $\xi$ is distributed uniformly on $[0,1]$ and is independent of the Wiener process.

\begin{theorem}
Let $(X_1^\e(t),X_2^\e(t))$ be the diffusion process in $\R^2$ defined above with initial conditions $X_1^\e(0) = x_1$ and $ X_2^\e(0) = x_2$. 
\begin{enumerate}[a)]
\item If $x_1 < 0$, then $(X_1^\e(t),X_2^\e(t))$ converges weakly as $\e \to 0$ on each finite time interval to the process $Y_t$ on $\Pi$ with $Y_0 = (x_1,x_2,4)$.
\item If $x_1 \geq 0$, then $(X_1^{\e,\d}(t),X_2^{\e,\d}(t))$ converges weakly as $\e \to 0$ on each finite time interval to the process $Y_t$ on $\Pi$ with a random initial condition $(x_1,x_2, \zeta)$, where $\zeta$ takes values $1$, $2$ and $3$ with probabilities $1+z_1-z_2$, $z_3-z_1$ and $1+z_2-z_3$, respectively. Note that the process $Y_t$ is independent of $\d$ so that the process $Y_t$ is actually the double limit of $(X_1^{\e,\d}(t), X_2^{\e,\d}(t))$ as first $\e \to 0$ and then $\d \to 0$.
\end{enumerate}
\end{theorem}
The proof of this theorem and more general results on homogenization can be found in \cite{CerraiPrep}. 

We briefly now consider perturbations of a completely integrable Hamiltonian system with more than one degree-of-freedom. Assume that one can introduce action-angle coordinates so that the perturbed system, after the time re-scaling has the form:
\begin{equation}\label{eq_book_hamil_per}
	\begin{cases}
	 \dot{I}_t^\e = \beta_1(I_t^\e,\varphi_t^\e), & I_0^\e = (I_1,I_2) \in \R^2,
	\\  \dot{\varphi}_t^\e = \frac{1}{\e} w(I_t^\e) + \beta_2(I_t^\e, \varphi_t^\e), & \varphi_0^\e = (\varphi_1,\varphi_2) \in \mathbb{T}^2,
	\end{cases}
\end{equation}
and the non-perturbed system is
\begin{equation}\label{eq_book_hamil_unper}
	 \dot{\bar{I}}_t = 0, \hs4 \dot{\bar{\varphi}}_t = w(\bar{I}_t).
\end{equation}
A two-dimensional torus corresponds to each $I_0 \in \R^2$. If $w_1(I_0)$ and $w_2(I_0)$ are rationally independent, just one invariant probability measure, the uniform distribution, is concentrated on this torus. But if $w_1(I_0)$ and $w_2(I_0)$ are rationally dependent, then the torus is covered by a family of periodic trajectories with one invariant measure concentrated on each such trajectory. In that case, $\Me$ for system $\eqref{eq_book_hamil_unper}$ can be parametrized by a modified two-dimensional domain $\mathcal{E}$ where the actions $I_1$ and $I_2$ are changing; each point of $\mathcal{E}$ for which $w_1(I)$ and $w_2(I)$ are rationally dependent should be replaced by a circle. Typically, the set 
\begin{equation*}
\mathcal{R} := \{(I_1,I_2) \in \mathcal{E}:w_1(I) \text{ and }w_2(I) \text{ are rationally dependent} \}, 
\end{equation*}
called the resonance set, is dense in $\mathcal{E}$, which makes the set $\Me$ rather complicated. However, if (the resonance set) $\mathcal{R}$ is small enough the problem can be regularized so that the evolution of $I_t^\e$ for $\e \ll 1$ can be described as a motion in $\mathcal{E}$. In the classical theory of dynamical systems, the regularization consists of considering the convergence in the Lebesgue measure of the space of initial conditions (see \cite{Anosov60}, \cite{Lochak88} and references therein). The effect of such a regularization is equivalent to the addition of appropriate small random perturbation to initial conditions. However, this type of regularization does not work if the collection of first integrals has critical points. Regularization by the addition of a noise in equations allows one to describe the evolution of $I_t^\e$ defined by $\eqref{eq_book_hamil_per}$ as a motion on the corresponding open book in more general situation (see \cite{Freidlin12}, Ch. 9). In this case certain gluing conditions on the binding of the open book should be imposed. These gluing conditions can lead to stochasticity of the long time behavior in pure deterministic systems (\cite{Freidlin12}, Ch. 9).

\section{Semi-flows in functional space}\label{sec_flows}
Consider a system of ordinary differential equations depending on a parameter $x \in \R^n$:
\begin{equation}\label{eq_flow_ode_nonper}
	\dot{u}_k(t,x) = f_k(x,u_1,...,u_m), 
	\hs4  u_k(0,x) = g_k(x),  \hs4  k \in \{1,...,m\}.
\end{equation}
We assume that functions $f_k(x,u)$ have bounded and continuous first and second derivatives and the functions $g_k(x)$ are bounded. We will look at $u(t,\cdot) = (u_1(t,\cdot)$, $u_2(t,\cdot),...,u_m(t,\cdot))$ as a semi-flow in the space of bounded measurable functions from $\R^n$ into $\R^m$. For any $k = 1,...,m$ let
\begin{equation*}
	L_kh(x) = \frac{1}{2} \sum_{i,j=1}^n a_{ij}^{(k)}(x) \frac{\partial^2 h}{\partial x_i \partial x_j} + \sum_{i=1}^n b_i^{(k)} (x) \frac{\partial h}{\partial x_i}.
\end{equation*}
As usual we assume that the coefficients $a_{ij}^{(k)} $ and $b_i^{(k)}$ are bounded together with their first and second derivatives, and that the matrices $(a_{ij}^{(k)})$ are non-negative definite. System $\eqref{eq_flow_ode_nonper}$ is our non-perturbed system. The perturbed system is then given by
\begin{equation}\label{eq_flow_ode_per}
	\dot{u}_k^\e(t,x) = f_k(x,u_1^\e, ...,u_m^\e) + \e L_k u_k^\e, \hs4 u_k^\e(0,x) = g_k(x). 
\end{equation}
Under mild additional conditions, system $\eqref{eq_flow_ode_per}$ has a unique solution $u^\e(t,x) = (u_1^\e,...,u_m^\e)$ that converges to $u(t,x)$ as $\e \downarrow 0$ on any finite time interval. Our goal is to describe behavior of $u^\e$ on time intervals growing together with $\e^{-1}$.

We consider the case where $m = 1$. Assume that the matrix $a^{(i)}(x) = a(x)$ is positive definite, $g(x) \geq 0$ and $f_1$ has the form $  f_1 = f(x,u) = u\cdot c(x,u)$ . Suppose first the non-linear term is of Fisher-Kolmogorov-Petrovskii-Piskunov (FKPP) type: that is, $c(x,u) > 0$ if $ u < 1$, $c(x,u) < 0$ if $u > 1$, and $c(x,0) = c(x) = \ds \max_{u \geq 0} c(x,u)$ (\cite{Fisher37}, \cite{Kolmogorov37}). Moreover, let $b_1(x) \equiv 0$ so that $L = L_1 = \frac{1}{2} \sum_{i,j=1}^n a_{ij}(x) \frac{\partial^2}{\partial x_i\partial x_j}$ and assume that $g$ is continuous in $\R^n$ except, maybe, at a finite set of smooth manifolds of dimension less than $n$. In this case the non-perturbed system has the form $\dot{u} = c(x,u) \cdot u$ while the perturbed system becomes 
\begin{equation*}
	\dot{u}^\e(t,x) = \e L u^\e + c(x,u) \cdot u,\hs4 u^\e(0,x) = g(x).
\end{equation*}
Let $G:= \mathrm{supp} \ g$. Due to our assumptions on $c(x,u)$, for $t$ large enough fixed and $\e > 0$ small, $u^\e(t,x)$ will be close to the step-function equal to $1$ inside the support $G$ of the initial function $g(x)$ and to $0$ outside of $G$. Such a step function is an equilibrium of the non-perturbed flow and the $\d$-measure concentrated on each such equilibrium is an ergodic invariant probability measure for the non-perturbed semi-flow. Moreover, the non-perturbed system has no other ergodic probability measures so that $\Me$ for this flow is the collection of $\d$-measures concentrated on step-functions with values $0$ and $1$. Each such step-function is defined by its support and a motion on $\Me$ can be described by the evolution of the support.

Now let $\tilde{u}^\e(t,x) = u^\e(\frac{t}{\sqrt{\e}},x).$ Then $\tilde{u}^\e(t,x)$ is the solution of the Cauchy problem
\begin{equation}\label{eq_flow_Cauchy}
\dot{\tilde{u}}^\e = \sqrt{\e} L \tilde{u}^\e + \frac{1}{\sqrt{\e}} c(x,\tilde{u}^\e) \tilde{u}^\e, \hs4 \tilde{u}^\e(0,x) = g(x).
\end{equation}
Consider the diffusion process $X_t^\e$ governed by the operator $\sqrt{\e} L$. This process satisfies
\begin{equation}\label{eq_flow_diffus}
X_t^\e - x = \e^{1/4} \int_0^t \sigma(X_s^\e) dW_s ,
\end{equation}
where $W_s$ is the Wiener process and $  \sigma(x) \sigma^*(x) = a(x)$. Using the Feynman-Kac formula, one can write down the following equation for $\tilde{u}^\e(t,x)$:
\begin{equation}\label{eq_flow_feynman}
	\tilde{u}^\e(t,x) = \E_x g(X_t^\e) \exp \Big[ \frac{1}{\sqrt{\e}} \int_0^t c(X_s^\e, \tilde{u}^\e(t-s,X_s^\e))ds\Big].
\end{equation}
It is easy to prove that system $\eqref{eq_flow_diffus}$, $\eqref{eq_flow_feynman}$ has a unique solution under our assumptions, and defined in this way function $\tilde{u}^\e$ solves the Cauchy problem $\eqref{eq_flow_Cauchy}$, at least if $g(x)$ is continuous or has simple discontinuities.

Since $c(x,u) \leq c(x,0) = c(x)$, $\eqref{eq_flow_feynman}$ implies that
\begin{equation}\label{eq_flow_bound}
\tilde{u}^\e(t,x)  \leq \sup_{x \in \R^n} |g(x)| \E_x \exp \Big\{ \frac{1}{\e} \int_0^t c(X_s^\e)dx\Big\}.
\end{equation}
Moreover, it follows from large deviation theory that
\begin{align}\label{eq_flow_devs}
\nonumber \lim_{\e \to 0} \sqrt{\e}  & \log \E_x \exp \Big[  \frac{1}{\sqrt{\e}} \int_0^t c(X_s^\e)ds\Big]
\\ & = \sup \Big\{ \int_0^t\Big[ c(\varphi_s) - \frac{1}{2} ( a^{-1}(\varphi_s) \dot{\varphi}_s \cdot \dot{\varphi_s} )\Big] ds : \varphi_0 = x , \varphi_t \in G\Big\} =: V_0(t,x).
\end{align}
First consider now the case where $c(x,u) = c(x) = c$ is constant. Then one can derive from $\eqref{eq_flow_bound}$ and $\eqref{eq_flow_devs}$ that 
\begin{equation*}
\lim_{\e \to 0} \tilde{u}^\e(t,x) = \begin{cases}
	0, &\text{if } \rho(x,G) > t \sqrt{2c},
	\\ 1, & \text{if } \rho(x,G) < t \sqrt{2c},
\end{cases}
\end{equation*}
where $\rho$ is the Riemannian metric corresponding to the quadratic form $ds^2 = (a^{-1}(x) dx\cdot dx)$. This means that the support $G_t$ of the $\ds \lim_{\e \to 0} \tilde{u}^\e(t,x)$ grows according to the Huygens principle with constant speed $\sqrt{2c}$ calculated in the metric $\rho(\cdot,\cdot)$ \cite{Freidlin79}.

When $c(x) \neq \mathrm{const}$, the situation is more sophisticated. Even if $x \in \R^1$, the function $t^*(x)$, defined by the equation $V_0(t^*(x),x) = 0$ is not necessarily monotone increasing with distance from $G$. This can lead to interesting new effects. For instance, the interface separating the areas where $\tilde{u}^\e(t,x)$ is close to $0$ and close to $1$ can have jumps \cite{Freidlin79}, \cite{Freidlin85}, \cite{Azencott13}. Such effects are interesting in models of biological evolution and in models for the propagation of infections. It is convenient to introduce a function
\begin{equation*}
	V_1(t,x) := \sup_{ \varphi_0 = x, \varphi_t \in G} \min_{0 \leq \theta \leq t} \int_0^\theta \Big[ c(\varphi_s) - \frac{1}{2} ( a^{-1}(\varphi_s) \dot{\varphi}_s \cdot \dot{\varphi}_s) \Big] ds.
\end{equation*}
One can check that $V_1(t,x) \leq 0$. Moreover, one can prove that $\ds \lim_{\e \to 0} \tilde{u}^\e(t,x) = 1$ inside the set $V_1(t,x) = 0$ and $\ds \lim_{\e \to 0} \tilde{u}^\e(t,x) = 0$ if $V(t,x) < 0$ (see \cite{Azencott13}, \cite{Freidlin91}, \cite{Freidlin96} and references therein). Most of these results were first proved by probabilistic methods and then later reproved (and sometimes improved) by various PDE based methods \cite{Evans89}, \cite{Barles90}.

One can also consider systems of reaction-diffusion equations under assumptions similar to the FKPP case, and describe the limiting motion on $\Me$ using limit theorems for large deviations \cite{Freidlin85}, \cite{Freidlin17}. In addition, long-time evolution of solutions of reaction-diffusion equations with bi-stable nonlinear term and slow spatial transport can be described as a motion on the corresponding simplex of invariant probability measures (see \cite{Freidlin85} and references there). Various other asymptotic problems concerning reaction-diffusion phenomena in incompressible media can be considered using the approach presented in this paper (compare with \cite{Freidlin02}, \cite{Freidlin95}, \cite{Freidlin18}, \cite{Freidlin08}).

Consider now another non-perturbed semi-flow:
\begin{equation}\label{eq_flow_integro}
	\dot{u}(t,x) = c(x,v_t) \cdot u(t,x),
	\hs4 v_t = \int_{\R^n} u(t,y)dy,
	\hs4 u(0,x) = g(x).
\end{equation}
Here the rate of change of $u(t,x)$ is not local, but rather depends on values of $u$ at other points in space. Such models arise naturally in many biological problems. We assume that $c(x,v)$ is Lipschitz-continuous and bounded, and that the initial function $g(x) \geq 0$ is continuous and has compact support $G_0$. Moreover, assume there exists a bounded Lipschitz-continuous function $\alpha(x) \geq 0$ such that $c(x,\alpha(x)) = 0$ while $c(x,v) > 0$ for $v < \alpha(x),$ and $c(x,v) < 0$ for $v > \alpha(x)$.

Our perturbed system, after an appropriate time change, has the form
\begin{equation}\label{eq_flow_nonloc_per}
	\dot{u}^\e = \frac{\sqrt{\e}}{2} \Delta u^\e + \frac{1}{\sqrt{\e}} c(x,v_t^\e) \cdot u^\e	,
	\hs4 v_t^\e = \int_{\R^n} u^\e(t,y)dy,
	\hs4 u^\e(0,x) = g(x).
\end{equation}
We will see that in the case of equation $\eqref{eq_flow_nonloc_per}$ the function $u^\e(t,x)$ is approximated for $0 < \e \ll 1$ not by a step function with a moving interface, but by a running spike which can move continuously or make jumps.

Let $X_t^\e$ be the diffusion process in $\R^n$ corresponding to the operator $\frac{\sqrt{\e}}{2} \Delta$. Using the Feynman-Kac formula, we obtain
\begin{align}\label{eq_flow_feynman_nonloc}
\nonumber & u^\e(t,x) = \E_x g(X_t^\e) \exp \Big[ \frac{1}{\sqrt{\e}} \int_0^t c(X_s^\e, v_{t-s}^\e)ds \Big],
\\ & v_t^\e = \int_{\R^n} u^\e(t,y)dy.
\end{align}
System $\eqref{eq_flow_feynman_nonloc}$ has a unique solution. This then implies that the solution of $\eqref{eq_flow_nonloc_per}$ exists and is unique.

Actually, it is convenient to consider the non-perturbed system $\eqref{eq_flow_integro}$ in a larger space, namely the space $B$ of finite Borel measures $m(\gamma)$ on $\R^n$ endowed with the weak topology. The time evolution is then defined by the equation 
\begin{equation}\label{eq_flow_measure_nonloc}
 	m_t(\gamma) = \int_\gamma \exp \Big[ \int_0^t c(x,v_s)ds \Big]m_0(dx), \hs4 v_t = m_t(\R^n).
\end{equation}
Equation $\eqref{eq_flow_measure_nonloc}$ has a unique solution for each finite measure $m_0(dx)$. This can be proved by successive approximations. Note that if $m_0(\gamma)$ has a density $g(x)$ with respect to the Lebesgue measure, then $\eqref{eq_flow_measure_nonloc}$ is equivalent to $\eqref{eq_flow_integro}$. 

It is easy to see that the measure $\d(x-y)\alpha(y)$ concentrated at $y \in \R^n$ is a stable equilibrium of the equation $\eqref{eq_flow_measure_nonloc}$. The unit measure (in the space $B$ of measures) is an invariant ergodic probability measure of semi-flow $\eqref{eq_flow_measure_nonloc}$ for each $y \in \R^n$.

To describe the limiting behavior of the solution $u^\e(t,x)$ of system $\eqref{eq_flow_feynman_nonloc}$, for each non-decreasing right continuous function $h: [0,T] \to \R$, we define
\begin{align*}
	V_h(t,x) & := \sup \Big\{ \int_0^t(c(\psi_s,h_s) - \frac{1}{2} |\dot{\psi}_s|^2) ds: 
	\\ & \psi:[0,T] \to \R^n \text{ absolutely continuous }, \psi_0 \in G_0, \psi_t = x
	\Big\}, 
\end{align*}
The function $V_h(t,x)$ is continuous for $t \in (0,T)$ and $x \in \R^n$, and it satisfies $\ds \lim_{t \to 0} V_h(t,x) = 0$ if $x \in G_0$ and $\ds \lim_{t \to 0} V_h(t,x) = -\infty$ if $x \notin G_0$. Let $\Lambda_T$ be the set of measurabe functions $\varphi:[0,T] \to \R^n$ such that $\alpha(\varphi_t)$ is non-decreasing, right continuous and $\ds \lim_{t \to 0} \alpha(\varphi_t)$ exists. A function $\varphi^* \in \Lambda_T$ is called a maximal solution of the equation
\begin{equation}\label{eq_flow_stat}
	V_{\alpha(\varphi)} (t,\varphi_t) = 0,
\end{equation}
whenever equation $\eqref{eq_flow_stat}$ is satisfied for $\varphi \equiv \varphi^*$ and for any $t \in (0,T]$ and $x \in \R^n$,
\begin{equation*}
	V_{\alpha(\varphi^*)}(t,x) \leq 0.
\end{equation*}
\begin{theorem}
	Assume that equation $\eqref{eq_flow_stat}$ has a unique maximal solution $\varphi^* \in \Lambda_T$. Then $\alpha(\varphi_t^*)$ is continuous everywhere on $[0,T]$ with the possible exception of (at most) a countable set $\Theta \subset [0,T]$. For $t \in [0,T] \setminus \Theta$, the solution $u^\e(t,x)$ of problem $\eqref{eq_flow_nonloc_per}$ converges weakly to the running spike $\alpha(\varphi_t^*)\delta(x-\varphi_t^*)$ as $\e \downarrow 0$.
\end{theorem}
The proof of this theorem is based on Theorem 1 from \cite{Freidlin87}. We note that in the interesting case where $n = 1$ and $ c(x,v) = \alpha(x) - v$, equation $\eqref{eq_flow_stat}$ can be solved in a sense, explicitly.

We also note that various symmetries can lead to the non-uniqueness of the maximal solution. In that case the limiting behavior of the solution of $\eqref{eq_flow_nonloc_per}$ may display branching and stochastic behavior, after an appropriate regularization \cite{Freidlin05}.

Consider now the following equation with a small delay:
\begin{equation}\label{eq_flow_delay_per}
\ddot{q}_t^\e = f(q_t^\e, q_{t-\e}^\e), \hs4 q_s^\e = \beta_s, \hs4 \dot{q_s}^\e = \dot{\beta}_s \text{ for }s \in [-\e,0].
\end{equation}
Here $\beta: [-\e, 0] \to \R$ is a continuously differentiable function, while the function $f(q,p)$ is assumed to be bounded together with its first and second derivatives. For small $\e \ll 1$, we can write
\begin{equation}\label{eq_flow_delay_approx}
	\ddot{q}_t^\e = f(q_t^\e,q_t^\e) - \e f_2'(q_t^\e,q_t^\e) \dot{q}_t^\e + o(\e),
\end{equation}
where $f_2'(q,q) = \frac{\partial f(q,p)}{\partial p} |_{p = q}$. Let $F'(q) = -f(q,q)$.

Now we introduce the following one degree-of-freedom oscillator with a friction term:
\begin{equation}\label{eq_flow_oscillator}
\ddot{\bar{q}}_t^\e = -F'(\bar{q}_t^\e) - \e f_2'(\bar{q}_t^e,\bar{q}_t^\e), \hs4 \bar{q}_0 = \beta_0, \hs4 \dot{\bar{q}}_0 = \beta_0'.
\end{equation}
Put $H(p,q)  = \frac{p^2}{2} + F(q)$. Let $\Gamma$ be the graph counting ergodic invariant probability measures of the non-perturbed system $\ddot{\bar{q}}_t = -F'(\bar{q}_t)$ corresponding to $\eqref{eq_flow_oscillator}$, and $Y :\R^2 \to \Gamma$ be the corresponding mapping (see Section \ref{sec_graph}). Then, $Y(\dot{\bar{q}}_{t/\e}^\e, \bar{q}_{t/\e}^\e)$ converges as $\e \downarrow 0$ to a certain motion $Y_t$ on $\Gamma$ (perhaps after regularization if $H(p,q)$ has saddle points). This follows from Ch. 8 of \cite{Freidlin12}. Combining this with $\eqref{eq_flow_delay_approx}$, one can prove that $Y(\dot{q}_{t/\e}^\e,q_{t/\e}^\e)$ converges to $Y_t$ as $\e \to 0$ on each finite time interval \cite{Freidlin05}. Note that the small delay can lead to a stochastic limiting motion $Y_t$ on $\Gamma $ in spite of the purely deterministic nature of equation $\eqref{eq_flow_delay_per}$. 

One can also consider perturbations of semi-flows defined by classical evolutionary PDE's. For example, one can add random perturbations to reaction-diffusion equations like the following:
\begin{equation}\label{eq_flow_reac}
\dot{u}^\e = \frac{1}{2} \Delta u^\e + F'(u^\e) + \e \mathcal{G}(t,x), \hs4 u^\e(0,x) =g(x), \hs4 x \in \R^1,
\end{equation}
where $\mathcal{G}(t,x)$ is a Gaussian space-time white noise. One can calculate the action functional and corresponding quasi-potential for equation $\eqref{eq_flow_reac}$ and similar RDE systems \cite{Faris82}, \cite{Freidlin88}. This allows one to describe the long-time behavior of solutions as $\e \downarrow 0$. If $x$ is multidimensional, one may have to assume that the noise has certain spatial regularity in $x$. This provides existence and uniqueness of a solution, but also makes expressions for the action and quasi-potential more complicated. Conversely, if the noise $\mathcal{G}(t,x)$ is close to the space-time white noise in a sense, one can obtain simpler expressions similar to the one-dimensional case \cite{Cerrai11}.

\section{Fast oscillating perturbations. Diffusion approximation.}\label{sec_oscil}
Consider a system in $\R^n$
\begin{equation}\label{eq_oscil_sys_per}
\dot{X}_t^\e = b(X_t^\e, h_{t/\e}), \hs4 X_0^\e = x \in \R^n,
\end{equation}
where the vector field $b:\R^n \times \R^m \to \R^n$ is assumed to be Lipschitz continuous and $h_t$ is a stationary stochastic process with sufficient mixing properties. Moreover, set $\bar{b}(x):=\E b(x,h_t) $ and let $R(s,t,x,y) = (R_{ij}(s,t,x,y))$ be the covariance matrix of the random vectors $b(x,h_s)$ and $b(y,h_t)$. Assume that there exists a function $\alpha:\R^+ \to \R^+$ such that $\ds \lim_{\tau \to \infty} \alpha(\tau) = 0$ and
\begin{equation*}
	\sup_{x,y \in \R^n, 1 \leq i,j \leq n} |R_{ij}(s,t,x,y)| \leq \alpha(|t-s|).
\end{equation*}
Under these conditions it is easy to check that $X_t^\e$ converges in probability as $\e \downarrow 0$ uniformly on each finite time interval to the solution $\bar{X}_t$ of the problem
\begin{equation}\label{eq_oscil_sys_unper}
	\dot{\bar{X}}_t = \bar{b}(X_t), \hs4 \bar{X}_0 = x.
\end{equation}
This law-of-large-numbers type result means that $X_t^\e$ can be considered as a fast oscillating perturbation of $\eqref{eq_oscil_sys_unper}$. According to our approach, to describe the long-time behavior of $\eqref{eq_oscil_sys_per}$, we should, first look at the simplex $\M$ of invariant probability measures of $\eqref{eq_oscil_sys_unper}$ and consider the projection of $X_t^\e$ on $\M$. Then we should describe the limiting evolution of this projection in an appropriate time scale.

The time scale of the evolution depends on the structure of the set $\Me$ of ergodic probability measures for system $\eqref{eq_oscil_sys_unper}$. If system $\eqref{eq_oscil_sys_unper}$ has a finite number of asymptotically stable ergodic probability measures, like in Section \ref{sec_finite}, the transition between different attractors occur in an exponential time scale. Similar to Section \ref{sec_finite}, we will have a hierarchy of cycles, metastable states, and results concerning exit problems; however, the action functional will be different here. For instance, if $h_t$ is a non-degenerate diffusion process on a compact manifold $M$ governed by an elliptic operator $L$, the action functional for the family $X_t^\e$ has the form
\begin{equation*}
	\frac{1}{\e} \int_0^t \ell(\varphi_s, \dot{\varphi}_s)ds,
\end{equation*}
where $\varphi_s$ is an absolutely continuous on $[0,T]$ function with values in $\R^n$ such that $\varphi_0 = X_0^\e = x$. The function $\ell(x,\alpha)$ is the Legendre transform in variable $\beta$ of the eigenvalue $\lambda(x,\beta)$ corresponding to a positive eigenfunction of the following problem on $M$:
\begin{equation*}
	L e(x,\beta,y) + (\beta \cdot \bar{b}(x)) e(x,\beta,y) = \lambda(x,\beta) e(x,\beta,y).
\end{equation*}
Such an eigenvalue exists, is simple and convex in $\beta \in \R^n$ (Ch. 7 of \cite{Freidlin12}). However, one should keep in mind that in this case, not all transitions between the attractors of system $\eqref{eq_oscil_sys_unper}$ may be possible.

Suppose now that equation $\eqref{eq_oscil_sys_unper}$ has a first integral $H(x)$. Moreover, assume that $H(x)$ is sufficiently smooth and satisfies $\ds \lim_{|x| \to \infty} H(x) = \infty$. Since $X_t^\e$ converges to $\bar{X}_t$, it follows that $H(X_t^\e) \to H(\bar{X}_t)$. Since $H(x)$ is a first integral for $\bar{X}_t$, we have $H(\bar{X}_t) \equiv H(x)$. Thus $\ds \lim_{\e \to 0} H(X_t^\e) = H(x)$ for any $t > 0$ independent of $\e$. To observe the evolution of $H(X_t^\e)$, we rescale the time by setting $\tilde{X}_t^\e := X_{t/\e}^\e$. This implies that
\begin{equation*}
\dot{\tilde{X}}_t^\e = \frac{1}{\e} b(\tilde{X}_t^\e, h_{t/\e^2}),
	\hs4 \tilde{X}_0^\e = x,
\end{equation*}
and
\begin{equation*}
	 H(\tilde{X}_t^\e) - H(x) = \frac{1}{\e} \int_0^t \nabla H(\tilde{X}_s^\e) \cdot (b(\tilde{X}_s^\e,h_{s/\e^2}) - \bar{b}(\tilde{X}_s^\e)) ds,
\end{equation*}
since $H(x)$ is a first integral so that $\nabla H(x) \cdot b(x) \equiv 0$. For each $x \in \R^n$, the quantity
\begin{equation*}
	\frac{1}{\e} \int_0^t \nabla H(x) \cdot (b(x,h_{s/\e^2}) - \bar{b}(x)) ds
\end{equation*}
converges as $\e \downarrow 0$ to a Gaussian random variable, provided $h_s$ has good enough mixing properties. Of course, the characteristics of this limiting random variable depend on $x$. Taking into account that $\tilde{X}_t^\e$ changes much slower than $h_{t/\e^2}$ and that $\tilde{X}_t^\e$ converges weakly to $\tilde{X}_t$ as $\e \to 0$, one can expect that if the dynamical system $\bar{X}_t$ has some ergodic properties on the level set $\{x: H(x) = y\}$, then the characteristics of the limit of $dH(X_t^\e)$ as $\e \to 0$ depend only on $H(X_t^\e)$. This means that the limiting process for $H(\tilde{X}_t^\e)$ as $\e \to 0$ will be the diffusion process
\begin{equation*}
	\dot{Y}_t = \sigma(Y_t) \dot{W}_t + B(Y_t), Y_0 = H(x),
\end{equation*}
where $W_t$ is a Wiener process. Thus the convergence of $H(\tilde{X}_t^\e)$ to a diffusion process is the result of mixing properties of the fast component $h_{t/\e^2}$ and of the ergodicity of the dynamical system $\bar{X}_t$ on the level sets. Exact assumptions on the mixing rate and rigorous results can be found in \cite{Borodin95}. In particular, explicit expressions for the diffusion and drift coefficients of the limiting process are calculated there.

The assumptions concerning the mixing are satisfied, if $h_t$ is a non-degenerate diffusion process on a compact manifold. For instance, they are satisfied if $h_t$ is the Wiener process on the $m$-dimensional torus, $\mathbb{T}^m$. In this case, mixing is exponentially fast so that one can take the process $h_t$ with any fixed initial point (rather than the stationary process).

The assumption concerning the ergodicity of $\bar{X}_t$ on the level sets, even if $n = 2$, is satisfied just in the case when $H(x)$ has one well. The set $\Me$ then can be parametrized by the values of $H$.

If $n = 2$ and $H(x)$ has several wells, as we have seen in Section $\ref{sec_graph}$, the set $\Me$ can be parametrized by a graph $\Gamma$. Inside the edges of the graph, one calculate the limiting diffusion process using the results mentioned above. But some of the vertices of $\Gamma$ can be accessible for the limiting diffusion process in a finite time. Certain gluing conditions should be added at these vertices. Sometimes, one can find these gluing conditions using the arguments mentioned before: Let system $\eqref{eq_oscil_sys_per}$ be a Hamiltonian system with one degree-of-freedom
\begin{equation*}
	\dot{X}_t^\e = \bar{\nabla } H(X_t^\e, h_{t/\e}), \hs4  X_0^\e = x \in \R^2,
\end{equation*}
where $h_t$ is the Wiener process on the unit circle $\mathbb{T}$. 

The pair $(X_t^\e,h_{t/\e})$ form a three-dimensional diffusion processs with a degenerate diffusion matrix. It is easy to check, for instance using the stationary forward Kolmogorov equation, that the Lebesgue measure on $\R^2 \times \mathbb{T}$ is invariant for the process $(X_t^\e, h_{t/\e})$ for each $\e > 0$. Let $\Gamma$ be the graph related to $\bar{H}(x) = \int_{\mathbb{T}} H(x,z)dz$ and $Y: \R^2 \to \Gamma$ be the corresponding projection (see Section \ref{sec_graph}). Since the Lebesgue measure in $\R^2 \times \mathbb{T}$ is invariant for $(X_t^\e, h_{t/\e})$, the projection $\lambda$ of the Lebesgue measure on $\Gamma$ induced by the mapping $Y$ is invariant for $Y(X_t^\e)$ (and also for $Y(\tilde{X}_t^\e)$) for each $\e >0$. This implies that the measure $\lambda$ is also invariant for the limit of $Y(\tilde{X}_t^\e)$ as $\e \downarrow 0$, if such a limit exists. Inside of each edge $I_k \subset \Gamma$, $Y(\tilde{X}_t^\e)$ converges to a diffusion process and the generator $L_k$ of this process can be calculated. It turns out that there exists just one diffusion process $Y_t$ on $\Gamma$ which is governed by the operators $L_k$ inside the edges and has the prescribed invariant measure $\lambda$. The gluing conditions for $Y_t$ can be expressed explicitly through the coefficients of operators $L_k$ and the measure $\lambda$.

To make these arguments rigorous, one should prove that the limiting process on $\Gamma$ is Markovian. To the best of my knowledge, it is not proved yet, although similar results for other problems are available.

Consider now the case when system $\eqref{eq_oscil_sys_unper}$ has several first integrals. To be specific, assume system $\eqref{eq_oscil_sys_per}$ has the form
\begin{equation*}
\dot{X}_t^\e = \bar{\nabla} H(X_t^\e, h_{t/\e}), \hs4 X_0^\e = x \in \R^{2n},
\end{equation*}
and the corresponding system $\eqref{eq_oscil_sys_unper}$ is a completely integrable Hamiltonian system with $n$ degrees of freedom of the form
\begin{equation}\label{eq_oscil_hamil_unper}
	\dot{\bar{X}}_t = \bar{\nabla}\bar{H}(\bar{X}_t), \hs4 \bar{X}_0 = x.
\end{equation}
Here $h_t$ is taken to be a Wiener process on $\mathbb{T}^m$. Let $H_1 = H(x),...,H_n(x)$ be the integrals of system $\eqref{eq_oscil_hamil_unper}$. One can consider the open book $\Pi$ homeomorphic to the set of connected components of all level sets $S_{z_1,...,z_n} = \{x \in \R^{2n}: H_1(x) = z_1,...,H_n(x) = z_n\}$. As was explained in Section \ref{sec_book}, because of the resonances, the set $\Me$ for the system $\eqref{eq_oscil_hamil_unper}$, in general, cannot be parametrized by the open book $\Pi$ if $n > 1$. But, if the resonance set is small enough, the limiting process for $(H_1(\tilde{X}_t^\e), ...,H_n(\tilde{X}_t^\e))$ as $\e \downarrow 0$ can be described as a diffusion process on $\Pi$, at least if the structure of $\Pi$ is simple enough. In particular, this can be done if $\Pi$ consists of one page. This question was considered in \cite{FreidlinArxiv}. It was shown there that the diffusion approximation for the $n$-dimensional process $(H_1(X_{t/\e}^\e),...,H_n(X_{t/\e}^\e))$ holds if the resonance set belongs to the union of a finite number of smooth surfaces of dimension less than $n$.

\section{Systems without finite invariant measures. Phantom dynamics.}\label{sec_phantom}

We now consider perturbations of systems that do not have finite invariant measures. A noisy perturbation in such a system can lead to the appearance of stable in a sense, unexpected attractors and other patterns. We call this phenomenon phantom dynamics. Similar effects can be observed in systems having unstable invariant manifolds and finite invariant measures concentrated on those manifolds.

To demonstrate this phenomenon, we consider a dynamical system in $\R^2$ of the form
\begin{equation}\label{eq_phantom_unper}
	\delta \dot{X}_t = f(X_t,Y_t), \hs4 \dot{Y}_t = X_t, \hs4 (X_0,Y_0) = (x,y) \in \R^2,
\end{equation}
where $f:\R^2 \to \R$ is assumed to be bounded together with its first and second derivatives and $\delta > 0$. Of course, this system is equivalent to the second order equation
\begin{equation}\label{eq_phantom_second_unper}
	\delta \ddot{Y}_t = f(\dot{Y}_t,Y_t), \hs4 (Y_0, \dot{Y}_0) = (x,y).
\end{equation}
\begin{figure}[h]
\centering
\includegraphics[scale = .8]{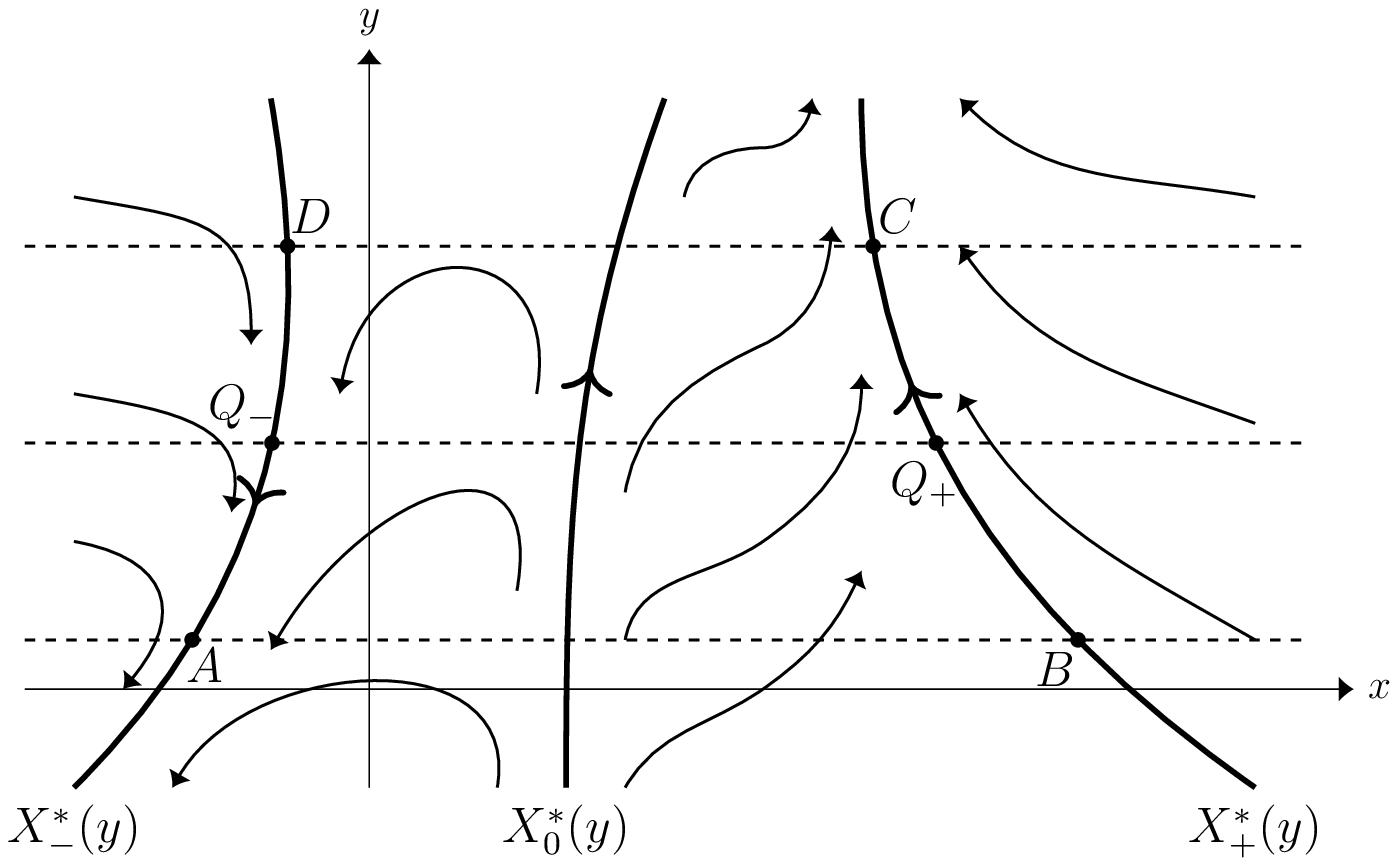}
\vspace*{-4mm}
\caption{}
\label{fig_phantom_vector}
\end{figure}

Suppose the set $\mathcal{E} = \{(x,y):f(x,y) = 0\}$ consists of three smooth curves $X_-^*(y) < 0 < X_0^*(y) < X_+^*(y),$ $y \in \R$, as shown in Figure \ref{fig_phantom_vector}. Moreover, suppose $f(x,y) > 0$ if the point $(x,y)$ is situation to the left of $X_-^*(y)$ or between $X_0^*(y)$ and $X_+^*(y)$, and $f(x,y) < 0$ for $(x,y)$ to the right of $X_+^*(y)$ or between $X_-^*(y)$ and $X_0^*(y)$. Then curves $X_-^*(y)$ and $X_+^*(y)$ attract the trajectories of system $\eqref{eq_phantom_unper}$, and the curve $X_0^*(y)$ repel the trajectories. It is clear that for any $\d > 0$, $Y_t$ tends to $+\infty$ or to $-\infty$ as $t \to \infty$, and the system has no finite invariant measures. Thus Assumption \ref{assume1} is not satisfied here.

Suppose now that we perturb the right hand side of $\eqref{eq_phantom_second_unper}$ by a small noise so that our system has the form 
\begin{equation}\label{eq_phantom_per}
	 \delta \dot{X}_t^\e = f(X_t^\e, Y_t^\e) + \sqrt{\e} \sigma(X_t^\e,Y_t^\e) \dot{W}_t, \hs4 \dot{Y}_t^\e = X_t^\e, \hs4 (X_0^\e,Y_0^\e) = (x,y),
\end{equation}
where the function $\sigma(x,y),(x,y) \in \R^2$ is assumed to be bounded, positive and Lipschitz-continuous.

The process $(X_t^\e,Y_t^\e)$ depends on two small parameters $\e$ and $\d$. Later, we will specify the relation between them, but one should note that if $0 < \d \ll 1$, then $x$-component in $\eqref{eq_phantom_per}$ changes with the rate $\d^{-1}$ while the $y$-component changes with the rate of order $1$. Consider now the first of equations $\eqref{eq_phantom_per}$ with a frozen variable $Y=y$:
\begin{equation}\label{eq_phantom_frozen}
	\dot{X}_t^{\e,y} = \frac{1}{\d} f(X_t^{\e,y},y) + \frac{\sqrt{\e}}{\d} \sigma(X_t^{\e,y},y) \dot{W}_t, \hs4 X_0^{\e,y} = x.
\end{equation}
According to our assumptions, the dynamical system $X_t^{0,y}$ has two asysmptotically stable equilibriums $X_-^*(y)$ and $X_+^*(y)$. They are separated by an unstable equilibrium at $X_0^*(y)$. If $\e \ll 1$ but still positive, then transitions between $X_-^*(y)$ and $X_+^*(y)$ are possible as large deviations of $X_t^{\e,y}$ from $X_t^{0,y}$. The quasi-potentials with respect to the equilibriums $X_+^*(y)$ and $X_-^*(y)$ for the process $X_t^{\e,y}$ (\cite{Freidlin12}) are given by
\begin{equation*}
	V_{\pm}(y) = 2\int_{X_0^*(y)}^{X_{\pm}^*(y)} \frac{f(x,y)dx}{\sigma^2(x,y)}.
\end{equation*}\
Let $\tau_+^{\e,\d}(y)$ $(\tau_-^{\e,\d}(y)$ ) be the (random) time of first transition from $X_+^*(y)$ to $X_-^*(y)$ (from $X_-^*(y)$ to $X_+^*(y)$) for the process defined by $\eqref{eq_phantom_frozen}$. Then, as it follows from \cite{Freidlin12}, $\tau_{\pm}^{\e,\d}(y) = \delta \tilde{\tau}_{\pm}^{\e,\d}$, and
\begin{equation}\label{eq_phantom_exittime}
	\lim_{\e \to 0} \e \log \tilde{\tau}_{\pm}^{\e,\d} (y) = V_{\pm}(y) \d.
\end{equation}
We used here that $\tilde{X}_t^{\e,y} = X_{t/\d}^{\e,y}$ satisfies the equation
\begin{equation*}
	\dot{\tilde{X}}_t^{\e,y} = f(\tilde{X}_t^{\e,y},y) + \sqrt{\frac{\e}{\d}} \sigma(\tilde{X}_t^{\e,y},y) \dot{W}_t.
\end{equation*}
\begin{figure}[h]

\centering
\includegraphics[scale = .9]{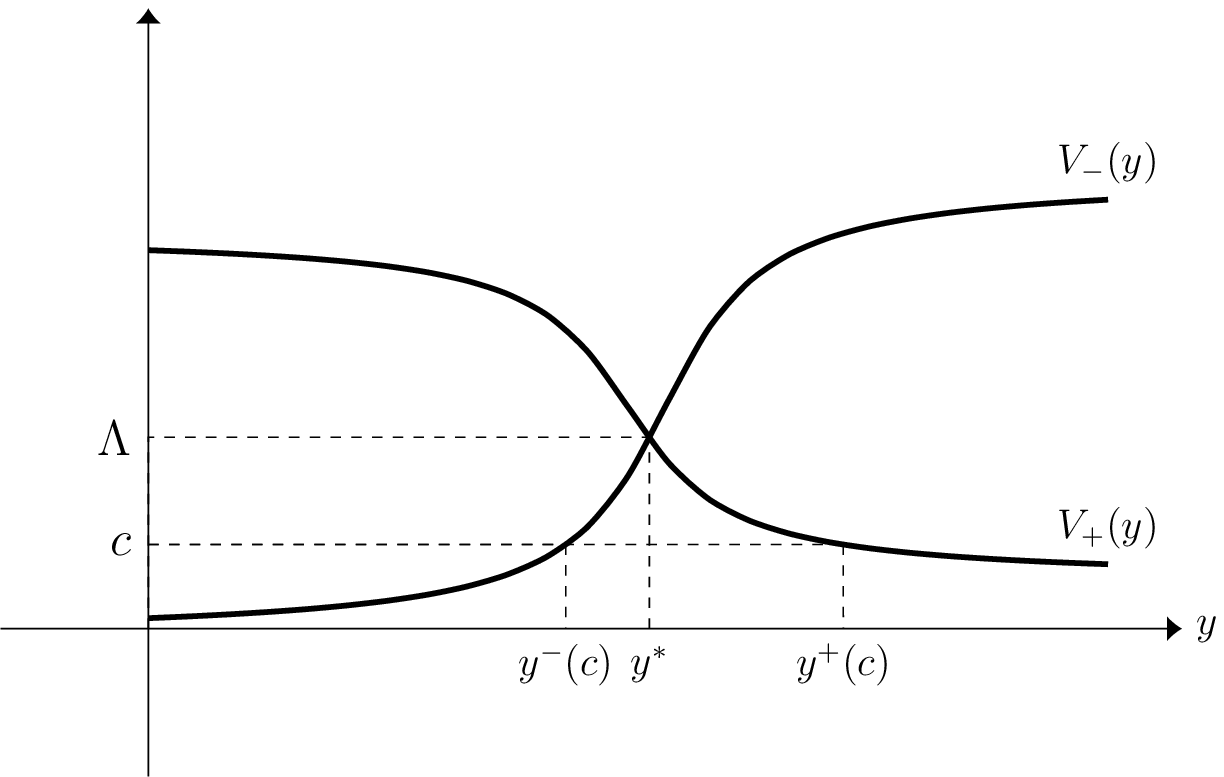}
\vspace*{-6mm}
\caption{}
\label{fig_phantom_quasi}
\end{figure}

Now, assume that the functions $V_{\pm}(y)$ are monotone, $ V_+'(y) < 0$, $V_-'(y) > 0$, $\bar{V}_{\pm} = \inf_{y \in \R} V_{\pm}(y) > 0$ and for some $y^*$ and $\Lambda$, $V_+(y^*) = V_-(y^*) = \Lambda $, as in Figure \ref{fig_phantom_quasi}. Moreover, assume that $\e$ and $\d$ are simultaneously sent to $0$ such that
\begin{equation}\label{eq_phantom_limit}
	\lim_{\e,\d \to 0} \frac{\e}{\d} = 0, \hs4 \liminf_{\e,\d \to 0} \frac{\e}{\d} \log\Big(\frac{1}{\d}\Big) > \Lambda.
\end{equation}
Equality $\eqref{eq_phantom_exittime}$ implies that if $\eqref{eq_phantom_limit}$ is satisfied, then for any $x \in \R$, $h > 0$ and $t > 0$ 
\begin{align}\label{eq_phantom_exit_asymp}
	\nonumber & \lim_{\e,\d \to 0} \tau_+^{\e,\d}(y)= 0, \hs4 \lim_{\e,\d \to 0} \mathbb{P}_x \{ |X_t^{\e,y}-X_-^*(y)| > h \} = 0 \text{ if } y > y^*,
	\\ &\lim_{\e,\d \to 0} \tau_-^{\e,\d}(y)= 0, \hs4 \lim_{\e,\d \to 0} \mathbb{P}_x \{ |X_t^{\e,y}-X_+^*(y)| > h \} = 0 \text{ if } y < y^*.
\end{align}
Taking into account that the $y$-component of the process $\eqref{eq_phantom_per}$ changes with the rate of order $1$, we can derive from $\eqref{eq_phantom_exit_asymp}$ that similar relations hold for the process $\eqref{eq_phantom_per}$:
\begin{align}\label{eq_phantom_exit_asymp_doub}
	\nonumber & \lim_{\e,\d \to 0} \mathbb{P}_{x,y} \{ |X_t^\e - X_-^*(y)| > h \} = 0, \ \text{ if }y > y^*,
	\\ &\lim_{\e,\d \to 0} \mathbb{P}_{x,y} \{ |X_t^\e - X_+^*(y)| > h \} = 0, \ \text{ if }y > y^*.
\end{align}
Since $f(X_-^*(y),y) > 0$ and $f(X_+^*(y),y) < 0$ for $y \in \R$, $\eqref{eq_phantom_exit_asymp_doub}$ implies that the process $(X_t^\e, Y_t^\e)$ for large enough $t$ (independent of $\e$ and $\d$) situates near points $\mathcal{Q}_- = (X_-^*(y^*), y^*)$ and $\mathcal{Q}_+ = (X_+^*(y),y)$ with probability close to $1$, provided $\eqref{eq_phantom_limit}$ holds (for more details see \cite{Freidlin01osc}).

Finally, we should calculate the limiting distribution $(P_-,P_+)$ between the points $\mathcal{Q}_-$ and $\mathcal{Q}_+$. Since $Y_t^\e$ should be close to $y^*$ for all large enough $t$, it follows that
\begin{equation*}
	P_-f(X_-^*(y^*),y^*) = P_+f(X_+^*(y),y^*).
\end{equation*}
Thus
\begin{align*}
	& P_+ = \frac{ f(X_-^*(y^*),y^*)}{f(X_-^*(y^*),y^*) + f(X_+^*(y^*),y^*)},
	\\ & P_- = \frac{ f(X_+^*(y^*),y^*)}{f(X_-^*(y^*),y^*) + f(X_+^*(y^*),y^*)}.
\end{align*}
We can summarize this result in the following theorem.
\begin{theorem}
	Let $(X_t^\e,Y_t^\e)$ be defined by equations $\eqref{eq_phantom_per}$, and the phase diagram of $(X_t^0,Y_t^0)$ be as shown in Figure \ref{fig_phantom_vector}. Let $y^*$ be the solution of the equation $V_+(y^*) = V_-(y^*)$ with $V_\pm(y)$ as shown in Figure \ref{fig_phantom_quasi}. Assume that conditions $\eqref{eq_phantom_limit}$ are satisfied. Then for any continuous bounded function $g:\R^2 \to \R$,
	\begin{equation*}
		\lim_{\e,\d \to 0, t \to \infty} \E_{x,y} g(X_t^\e, Y_t^\e) = P_-g(X_-^*(y^*),y^*) + P_+g(X_+^*(y^*),y^*).
	\end{equation*}
\end{theorem}
The noise induced a probability measure concentrated on the set consisting of two points $\mathcal{Q}_-$ and $\mathcal{Q}_+$ such that the trajectory is attracted to this measure. The value $y^*$ is effectively a stable equilibrium for equation $\eqref{eq_phantom_second_unper}$. If equation $V_+(y) = V_-(y)$ has many solutions, the system can have several such stable equilibriums. In larger time scales, transitions between these equilibriums due to large deviations are possible. One can give conditions leading in the limit to stable oscillations or to other patterns which are not available in the non-perturbed system. Systems with many degrees-of-freedom and other types of noise can be considered as well (compare with \cite{Freidlin01osc}, \cite{Freidlin01stoc}).

\end{document}